\documentclass[11pt,a4paper,twoside,enabledeprecatedfontcommands,final]{scrartcl}
\usepackage{a4wide}
\usepackage{amsfonts}
\usepackage{amsmath}
\usepackage{amssymb}
\usepackage{amsthm}
\usepackage[numbers]{natbib}
\usepackage{booktabs}
\usepackage{algorithm}
\usepackage{algorithmic}
\usepackage{pgfplots}
\pgfplotsset{compat=1.8}
\usepackage{graphicx}
\usepackage{color}
\usepackage{colortbl}   
\usepackage{showkeys}
\usepackage{scrpage2}
\usepackage{multirow}
\usepackage{boxedminipage}
\usepackage{enumerate}
\usepackage{subfig}
\usepackage{textcomp}

\usepackage{hyperref}

\usepackage{todonotes}
\usepackage{placeins}

\usepackage{siunitx}
\newcommand{\N}{\ensuremath{\mathbb{N}}}

\newcommand{\T}{\ensuremath{\mathbb{T}}}

\newcommand{\Z}{\ensuremath{\mathbb{Z}}}

\newcommand{\R}{\ensuremath{\mathbb{R}}}
\newcommand{\C}{\ensuremath{\mathbb{C}}}

\newcommand{\ii}{\textnormal{i}}
\newcommand{\e}{\textnormal{e}}

\newcommand{\ceil}[1]{\left\lceil#1\right\rceil}

\newcommand{\zb}[1]{\ensuremath{\boldsymbol{#1}}}

\renewcommand{\ln}{\mathrm{ln\,}}

\newcommand{\boldx}{{\ensuremath{\boldsymbol{x}}}}
\newcommand{\boldk}{{\ensuremath{\boldsymbol{k}}}}
\newcommand{\boldh}{{\ensuremath{\boldsymbol{h}}}}

\newcommand{\boldz}{{\ensuremath{\boldsymbol{z}}}}

\newcommand{\boldi}{{\ensuremath{\boldsymbol{i}}}}

\newcommand{\boldone}{{\ensuremath{\boldsymbol{1}}}}

\definecolor{darkgreen}{rgb}{0.0,0.5,0.0}

\DeclareMathOperator*{\esssup}{ess\,sup}

\newtheorem{theorem}{Theorem}[section]
\newtheorem{lemma}[theorem]{Lemma}
\newtheorem{remark}[theorem]{Remark}
\newtheorem{generalisation}[theorem]{Generalisation}
\newtheorem{definition}[theorem]{Definition}
\newtheorem{example}[theorem]{Example}
\newtheorem{corollary}[theorem]{Corollary}
\newtheorem{proposition}[theorem]{Proposition}

\newenvironment{Theorem}{\goodbreak \begin{theorem}\sl}{\end{theorem}}
\newenvironment{Lemma}{\goodbreak \begin{lemma}\sl}{\end{lemma}}
\newenvironment{Remark}{\goodbreak \begin{remark}\rm}{\bend\end{remark}}

\newenvironment{Example}{\goodbreak \begin{example}\rm}{\bend\end{example}}

\newenvironment{Corollary}{\goodbreak \begin{corollary}\sl}{\end{corollary}}

\makeatletter
\def\imod#1{\allowbreak\mkern10mu({\operator@font mod}\,\,#1)}
\makeatother

\numberwithin{equation}{section}
\numberwithin{table}{section}
\numberwithin{figure}{section}

\newcommand{\bend}{\hspace*{0ex} \hfill \hbox{\vrule height
    1.5ex\vbox{\hrule width 1.4ex \vskip 1.4ex\hrule  width 1.4ex}\vrule
    height 1.5ex}}

\long\def\symbolfootnote[#1]#2{\begingroup
\def\thefootnote{\fnsymbol{footnote}}\footnote[#1]{#2}\endgroup}

\newcommand{\OO}[1]{\mathcal{O}\left(#1\right)}

\renewcommand{\mathbf}[1]{\ensuremath{\boldsymbol{#1}}}
\renewcommand{\textbf}[1]{{\ensuremath{\boldsymbol{#1}}}}
\renewcommand{\thefootnote}{\fnsymbol{footnote}}

\title{High-dimensional sparse FFT based on sampling along multiple rank-1 lattices}
\date{}
\author{Lutz K\"ammerer\footnotemark[1] \and Daniel Potts\footnotemark[2] \and 
        Toni Volkmer\footnotemark[3]}
        
        \hypersetup{
          pdfcreator = {pdflatex},
          plainpages=false,
          pdfstartview=Fit,       
          pdfpagemode=UseOutlines, 
          bookmarksnumbered=true, 
          bookmarksopen=false,     
          bookmarksopenlevel=0,   
          colorlinks=true,       
          linkcolor=black,
          citecolor=black,
          urlcolor=black}

\begin{document}

\maketitle

\footnotetext[1]{
  Chemnitz University of Technology, Faculty of Mathematics, 09107 Chemnitz, Germany\\
  lutz.kaemmerer@mathematik.tu-chemnitz.de, Phone:+49-371-531-37728, Fax:+49-371-531-837728
}
\footnotetext[2]{
  Chemnitz University of Technology, Faculty of Mathematics, 09107 Chemnitz, Germany\\
  potts@mathematik.tu-chemnitz.de, Phone:+49-371-531-32150, Fax:+49-371-531-832150
}
\footnotetext[3]{
  Chemnitz University of Technology, Faculty of Mathematics, 09107 Chemnitz, Germany\\
  toni.volkmer@mathematik.tu-chemnitz.de, Phone:+49-371-531-39999, Fax:+49-371-531-839999
}

\begin{abstract}

The reconstruction of high-dimensional sparse signals is a challenging task in a wide
range of applications.
In order to deal with high-dimensional problems,
efficient sparse fast Fourier transform algorithms
are essential tools.
The second and third authors have recently proposed a dimension-incremental approach,
which only scales almost linear in the number of required sampling values and 
almost quadratic in the arithmetic complexity with respect to the spatial dimension $d$.
Using reconstructing rank-1 lattices as sampling scheme, the method
showed reliable reconstruction results in numerical tests but
suffers from relatively large numbers of samples and arithmetic operations.
Combining the preferable properties of reconstructing rank-1 lattices
with small sample and arithmetic complexities, the first author
developed the concept of multiple rank-1 lattices.
In this paper, both concepts  --- dimension-incremental reconstruction and
multiple rank-1 lattices --- are coupled, which yields a distinctly improved
high-dimensional sparse fast Fourier transform.
Moreover, the resulting algorithm is analyzed in detail with respect to
success probability, number of required samples, and arithmetic complexity.
In comparison to single rank-1 lattices, the utilization of multiple rank-1 lattices
results in a reduction in the complexities by an almost linear factor with respect to the sparsity.
Various numerical tests confirm the theoretical results, the high performance,
and the reliability of the proposed method.

\medskip

\medskip

\end{abstract}
\newpage

\section{Introduction}

Equally weighted cubature formulas, called quasi Monte-Carlo rules, 
are substantial tools in the field of numerical integration.
One kind of such cubature formulas are so called lattice rules that
allow for the very efficient numerical treatment of high-dimensional integrals. In particular rank-1 lattices are very well investigated in this field of mathematics, cf. \cite{SlJo94,JoKuSl13} for more details on this topic.

In a natural way, different researchers used lattice rules in order to
approximate integrals that compute Fourier coefficients of functions, and thus, reconstruct periodic signals.
This reconstruction can be done by interpolation, cf. e.g. \cite{MuSo12, Be13}, as well as approximation, cf. \cite{Tem86, KuSlWo06, KuSlWo08, KuWaWo09, ByKaUlVo16}.
In particular, for highly interesting types of functions, e.g. functions of dominating mixed smoothness, suitable approximations use lattices that necessarily require a relatively large number of sampling nodes in comparison to the achievable worst case approximation
errors, cf. \cite{ByKaUlVo16}. A recent idea \cite{Kae16} to prevent the usage of a single huge rank-1 lattice is to sample along multiple rank-1 lattices in order to reconstruct trigonometric polynomials or
approximate functions. With high probability, this ansatz significantly reduces the number of required sampling values, cf. \cite{Kae17}.
The crucial innovation is the strategy that determines rank-1 lattices at random, which is also used in a more recent paper for numerical integration, cf. \cite{KrKuNuUl17}.

In this paper, we deal with the reconstruction and approximation of multivariate periodic functions based on samples without knowing the locations of the non-zero or approximately largest Fourier coefficients. A method which allows for efficiently performing this task will be called \textit{multivariate sparse fast Fourier transform} in the following.
One aim is to exactly reconstruct the Fourier coefficients $\hat{p}_\boldk$, $\boldk\in I$,
of an arbitrarily chosen trigonometric polynomial
\begin{equation}\label{poly}
p(\boldx):=\sum_{\boldk\in I} \hat{p}_\boldk \, \mathrm{e}^{2\pi\mathrm{i}\boldk\cdot\boldx}
\end{equation}
with frequencies supported on an index set $I\subset\Z^d$ of finite cardinality, i.e. $|I|<\infty$, from sampling values of $p$,
where the frequency index set $I$ is unknown. The Fourier coefficients $\hat{p}_\boldk\in\C$
of the trigonometric polynomial $p$
are formally given by
the Fourier transform of $p$,
\begin{equation}\label{polyfk}
 \hat{p}_\boldk := \int_{\T^d} p(\boldx) \, \mathrm{e}^{-2\pi\mathrm{i}\boldk\cdot\boldx} \mathrm{d}\boldx, \quad \boldk\in I.
\end{equation}
Alternatively, in case of a general function $f\in L_1(\T^d)\cap\mathcal{C}(\T^d)$, we aim to determine an approximation of the (roughly) largest Fourier coefficients
$$
 \hat{f}_\boldk := \int_{\T^d} f(\boldx) \, \mathrm{e}^{-2\pi\mathrm{i}\boldk\cdot\boldx} \mathrm{d}\boldx, \quad \boldk\in I,
$$
from sampling values of $f$.
The critical part for both tasks is the determination of the unknown frequency index set $I$. Using a tensor-product approach and computing a full-dimensional fast Fourier transform (FFT) is not feasible in practice already for medium dimensions (like $d=5$) due to the curse of dimensionality.

Various methods for multivariate sparse FFTs exist which are based on different ideas. In general these methods require a search domain $\Gamma\subset\Z^d$ of possible frequencies as an input parameter, which is often chosen as a full grid $\hat{G}_N^d:=\{-N,-N+1,\ldots,N\}^d$ of expansion $N\in\N$, $|\hat{G}_N^d|=(2N+1)^d$.
Usually, an additional input parameter is
the number of approximately largest Fourier coefficients $\hat{f}_\boldk$ to determine
or, in the case of the reconstruction of a multivariate trigonometric polynomial $p$,
an upper bound on the cardinality $|\mathrm{supp}\,\hat{p}|$ of the support in frequency domain $$\mathrm{supp}\,\hat{p}:=\{\boldk\in\Z^d\colon \hat{p}_\boldk\neq 0\}.$$
Next, we briefly mention some of the existing methods.

One approach, which requires a relatively small amount of samples,
is applying random sampling in compressed sensing \cite{Do06,Ca06,FoRa13},
but the number of required arithmetic operations contains the cardinality~$|\Gamma|$ of the search domain~$\Gamma$ as
a linear factor,
see e.g.~\cite{FoRa13,KaKuMePoVo14,KuRa06}, and therefore typically suffers heavily from the curse of dimensionality.

In~\cite{InKa14}, a multivariate sparse FFT method was discussed which uses randomized sub-sampling and filters. The method is based on the one-dimensional versions from~\cite{HaInKaPr12,HaInKaPr12a,InKaPr14}. As search domain~$\Gamma$, a full grid~$\hat{G}_N^d$ is used and
the number of required samples is $\mathcal{O}(\vert\mathrm{supp}\,\hat{p}\vert\,\log N)$
for constant dimension~$d$ and the arithmetic complexity is~$\mathcal{O}(N^d \log^{\mathcal{O}(1)} N)$.
We remark that the sample complexity~$\mathcal{O}(\vert\mathrm{supp}\,\hat{p}\vert\,\log N)$ 
might contain a factor of~$d^{\mathcal{O}(d)}$, cf.~\cite[Section~IV]{InKa14}.
This means that both complexities suffer from the curse of dimensionality.

In~\cite{Iw13}, a deterministic multivariate sparse FFT algorithm was presented, which uses the Chinese Remainder Theorem and
which requires $\mathcal{O}(d^4 \, \vert\mathrm{supp}\,\hat{p}\vert^2\log^4(dN))$ samples and arithmetic operations.
This means there is neither exponential/super-exponential dependency on the dimension $d\in\N$
nor a dependency on a failure probability in the asymptotics of the number of samples and arithmetic operations for this method.
Besides this deterministic algorithm, there also exists a randomized version which only requires
$\mathcal{O}(d^4 \, \vert\mathrm{supp}\,\hat{p}\vert\log^4(dN))$ samples and arithmetic operations with high probability.

In this work, we focus on a dimension-incremental method for a high-dimensional sparse FFT published in~\cite{PoVo14}.
The method adaptively constructs the frequency index set $I$ belonging to the non-zero or approximately largest Fourier coefficients in a dimension-incremental way based on sampling values.
One preferable property of this strategy is that the sampling complexity is linear in the dimension~$d$ up to logarithmic factors.
Another advantage is the applicability in practice for relative high sparsities even in the case where the sampling values are perturbed by noise.
Moreover, the dimension-incremental strategy provides a straightforward implementation for the (approximate) recovery of multivariate periodic signals.
Known properties of the function under consideration, such as smoothness conditions as well as (possibly very large) supersets of the locations of the non-zero or approximately largest Fourier coefficients,
can be directly exploited by restricting the search domain~$\Gamma$ of considered frequencies.
In general, the dimension-incremental approach is a probabilistic method, which repeatedly performs the sampling using $r\in\N$ many detection iterations.
Theoretical bounds on $r$ are currently unknown and will be addressed in this paper.

One crucial challenge of the dimension-incremental method is the construction of spatial discretizations
for high-dimensional sparse trigonometric polynomials that allow for an efficient
reconstruction. 
In~\cite{PoVo14}, the authors mainly used single rank-1 lattices, which worked very reliably in various numerical experiments.
Besides this reliability, the main advantage of using single rank-1 lattices as spatial discretizations is the
available fast Fourier transform for these sampling sets. However, the bottlenecks of this approach are
the relatively expensive construction of reconstructing rank-1 lattices and, in addition, the relatively
large number of sampling nodes.
In detail,
for sparsity $s:=|I|$ or when determining $s$ many approximate Fourier coefficients, the method requires
$\mathcal{O}(d\,r^3s^2N)$ samples and
$\mathcal{O}(d\,r^3s^3 + d\,r^3s^2N\log(r\,s\,N))$ arithmetic operations in the case $\sqrt{N}\lesssim s\lesssim N^d$, when restricting the search space $\Gamma$ of considered frequency candidates to a full grid $\hat{G}_N^d$ of expansion $N\in\N$.
In order to reduce these complexities, the paper \cite{PoVo14} also numerically investigated alternative sampling sets
that allow for efficient reconstructions of sparse signals. The authors observed in the numerical experiments,
that the additionally investigated
types of sampling sets do not offer the reliability and stability as obtained for single rank-1 lattices.

Again, we emphasize that the dimension-incremental idea from \cite{PoVo14} should work with any sampling method which reliably computes Fourier coefficients in a fast way.
Accordingly, the dimension-incremental sparse FFT may be improved by providing sampling sets that combine four
preferable properties:
\begin{itemize}
\item low oversampling factors,
\item stability and thus reliability,
\item efficient construction methods,
\item fast Fourier transform algorithms.
\end{itemize}

Exactly this combination of objectives motivated the development of the concept of reconstructing multiple rank-1 lattices in \cite{Kae16,Kae17}.
Therein, a reconstruction approach for multivariate trigonometric polynomials $p$ with known frequency index sets $I$ based on samples along multiple rank-1 lattices is proposed.
Instead of using a single rank-1 lattice for the reconstruction, which requires $\mathcal{O}(d\,|I|^3)$ arithmetic operations for the construction of a suitable reconstructing rank-1 lattice that consists of up to $\mathcal{O}(|I|^2)$ samples, a union of cleverly chosen rank-1 lattices is utilized, where the latter is named reconstructing multiple rank-1 lattice.
Such a reconstructing multiple rank-1 lattice consists of $\mathcal{O}(|I|\log |I|)$ nodes with high probability under mild assumptions and, in this case, the construction only requires $\mathcal{O}(|I|(d+\log|I|)\log |I|)$ arithmetic operations, cf. \cite[Corollary~3.7 and Algorithm~4]{Kae17}.

\begin{sloppypar}
As mentioned above, we include the multiple rank-1 lattice idea into the dimension-incremental reconstruction approach in this work. We obtain a method, see Algorithm~\ref{algo:sfft_mlfft}, which requires
$\mathcal{O}(dr^2sN \log (rsN))$
samples and
$\mathcal{O}\big(d^2r^2sN \, \log^2(rsN) \big)$
arithmetic operations.

One of the main contributions of this work is the detailed investigation of
the overall success probability of the dimension-incremental reconstruction method.
For the first time, we show bounds on the number of detection iterations~$r$ when using single and multiple rank-1 lattices as sampling sets.
When reconstructing a trigonometric polynomial $p$ using Algorithm~\ref{algo:sfft_mlfft}, we prove in Theorem~\ref{thm:r_complexity} that under mild assumptions and
when choosing the number of detection iterations $r:=\lceil 2|\operatorname{supp}\hat{p}|(\log3+\log d+\log |\operatorname{supp}\hat{p}|-\log\varepsilon) \rceil$ for a given failure probability $\varepsilon\in (0,1)$ ensures that the proposed algorithm successfully detects all non-zero Fourier coefficients $\hat{p}_\boldk\neq 0$
with probability at least $1-\varepsilon$. In this case, we require
$\OO{d\,|\operatorname{supp}\hat{p}|^2N\log^2 (|\operatorname{supp}\hat{p}|)\,|\log\varepsilon|}$ samples and $\OO{d^2 \,|\operatorname{supp}\hat{p}|^2N \,\log^3 (|\operatorname{supp}\hat{p}|)\,|\log\varepsilon|}$ arithmetic operations. In our numerical tests, 
we did not require the linear factor $|\operatorname{supp}\hat{p}|$ in the number of detection  iterations $r$, but
we were able to choose this parameter as a constant $\leq 10$.
It is an open question if a smaller bound on $r$ can be theoretically shown under mild assumptions on the Fourier coefficients $\hat{p}_\boldk$. The key result which would need to be improved is currently based on a Markov-type inequality in Lemma~\ref{lem:prob_g_ineq}, see also the discussion in Section~\ref{sec:theo:discussion}.
Still, the results in Lemma~\ref{lem:prob_g_ineq}, Corollary~\ref{cor:prob_AtC} and Corollary~\ref{thm:step_t_single_frequency} imply a distinct improvement compared to the previous ones in \cite[Lemma~2.4 and Theorem~2.5]{PoVo14}, which would yield a factor of $|\operatorname{supp}\hat{p}|^2$ in $r$ when applying the theoretical considerations in Section~\ref{sec:theoretical:prob} on \cite[Theorem~2.5]{PoVo14}.
\end{sloppypar}

\begin{sloppypar}
If the Fourier coefficients $\hat{p}_\boldk$ fulfill certain stronger properties, we obtain distinctly smaller sample and arithmetic complexities by slightly modifying the algorithms, see also Remark~\ref{rem:deterministic}.
As discussed in Section~\ref{sec:theo:deterministic}, if we know that the real parts of all non-zero $\hat{p}_\boldk$ have identical sign $\neq 0$ or that the imaginary parts of all non-zero $\hat{p}_\boldk$ have identical sign $\neq 0$,
then we can use a modified version of Algorithm~\ref{algo:sfft_mlfft} with $r:=1$ detection iteration.
The resulting algorithm successfully detects all non-zero Fourier coefficients $\hat{p}_\boldk\neq 0$ with probability at least $1-\varepsilon$,
requiring $\mathcal{O}\big(d\,|\operatorname{supp}\hat{p}|N \log (|\operatorname{supp}\hat{p}|\,N)\big)$ samples and $\mathcal{O}\big(d^2 \, |\operatorname{supp}\hat{p}|N \log^2(|\operatorname{supp}\hat{p}|\,N) \,(\log d - \log\varepsilon)\big)$ arithmetic operations, cf. Section~\ref{sec:theo:deterministic}.
\end{sloppypar}

The remaining parts of this paper are structured as follows. In Section~\ref{sec:pre}, we briefly review important aspects of reconstructing multiple rank-1 lattices and dimension-incremental reconstruction. The new method is discussed in Section~\ref{sec:method}. In Section~\ref{sec:theoretical}, we develop bounds on the number of detection iterations $r$ and on the success probability $1-\varepsilon$ for the reconstruction of multivariate trigonometric polynomials $p$ when using reconstructing single or multiple rank-1 lattices as sampling sets.
Numerical tests performed in Section~\ref{sec:numerics} demonstrate the high performance of the proposed method for high-dimensional trigonometric polynomials~$p$ as well as for a 10-dimensional test function~$f$ which is non-sparse in frequency domain. Finally, in Section~\ref{sec:conc}, we conclude the results of this paper.

\section{Prerequisites}
\label{sec:pre}

For our method, we combine two main ingredients. One is using so-called reconstructing multiple rank-1 lattices for known frequency index sets~$I$ from~\cite{Kae16,Kae17}. The second one is the dimension-incremental reconstruction approach for unknown frequency index sets~$I$ as presented in~\cite{PoVo14}.

\subsection{Reconstructing multiple rank-1 lattices}
\label{sec:pre:multiple_r1l}

First, we start with definitions from \cite{Kae17} using slightly adapted symbols in this work.

The sampling sets $\mathcal{X}$ constructed and used by the method presented in Section~\ref{sec:method} are based on so-called
rank-1 lattices
$$
\Lambda(\boldz,M):=\{j\boldz/M\bmod \boldone\colon j=0,\ldots,M-1 \} \subset \T^d,
$$
where
$\boldz\in\Z^d$ and $M\in\N$ are called generating vector and lattice size of $\Lambda(\boldz,M)$, respectively.
For an arbitrary multivariate trigonometric polynomial~$p$, see \eqref{poly}, with frequencies~$\boldk$ supported on an index set $I\subset\Z^d$, $|I|<\infty$, we can reconstruct all the Fourier coefficients $\hat{p}_\boldk$, $\boldk\in I$ given in \eqref{polyfk},
from samples along a rank-1 lattice $\mathcal{X}:=\Lambda(\boldz,M)$ if the Fourier matrix $$
\boldsymbol{A}(\mathcal{X},I):=\left( 
\mathrm{e}^{2\pi\mathrm{i}\boldk\cdot\boldx} \right)_{\boldx\in \mathcal{X}, \boldk \in I}
$$ has full column rank.
This is the case if and only if $\Lambda(\boldz,M)$ is a \textit{reconstructing rank-1 lattice} for $I$, i.e., the \textit{reconstruction property}
$$
 \boldk\cdot\boldz\not\equiv \boldsymbol{k'}\cdot\boldz\imod{M} \text{ for all } \boldk,\boldsymbol{k'}\in I, \boldk\neq\boldsymbol{k'},
$$
is fulfilled. Such a reconstructing rank-1 lattice $\Lambda(\boldz,M)$ can be easily constructed using a simple component-by-component construction method, cf. \cite{Kae2013}.
However, the construction method has rather high computational costs and may require $\mathcal{O}(d|I|^3)$ arithmetic operations.
Moreover, under mild assumptions, the lattice size~$M$ is bounded by $|I|\leq M \leq |I|^2$, where this number~$M$ tends more to the upper bound for many interesting structures of frequency index sets~$I$.

Recently, in \cite{Kae16,Kae17}, a modified approach was presented, which allows for drastically reducing the number of samples. This approach uses rank-1 lattices $\Lambda(\boldz,M)$ as building blocks and builds sampling sets based on multiple instances. The corresponding sampling sets are called multiple rank-1 lattices and they can be constructed by simple and efficient randomized construction algorithms.
A multiple rank-1 lattice is the union of $L\in\N$ many rank-1 lattices,
$$
\Lambda=\Lambda(\boldz_1,M_1,\ldots,\boldz_L,M_L):=\bigcup_{\ell=1}^L \Lambda(\boldz_\ell,M_\ell),
$$
and consists of
$$
M:=|\Lambda(\boldz_1,M_1,\ldots,\boldz_L,M_L)|
\leq
1-L+\sum_{\ell=1}^{L}M_\ell
$$
many distinct nodes.
If $\Lambda=\Lambda(\boldz_1,M_1,\ldots,\boldz_L,M_L)$ allows for the reconstruction of all multivariate trigonometric polynomials $p$ with frequencies supported on a frequency index set $I$,
it will be called \textit{reconstructing multiple rank-1 lattice} for $I$.

In simplified terms, the basic idea is that each of the rank-1 lattices $\Lambda(\boldz_\ell,M_\ell)$, $\ell=1,\ldots,L$, should be a reconstructing one for some index set $I_\ell\subset I$
and that $\bigcup_{\ell=1}^L I_\ell = I$. We remark that this condition is not sufficient in general and we require an additional property.

One construction method is Algorithm~\ref{alg:construct_mr1l_I_distinct_primes:mod}, which searches for a reconstructing multiple rank-1 lattice $\Lambda$ for a given index set $I$
such that the properties
\begin{equation}\label{equ:cond:reco_ml_mlfft2}
 \boldk\cdot\boldz_\ell\not\equiv \boldsymbol{k'}\cdot\boldz_\ell\imod{M_\ell} \text{ for all } \boldk\in I_\ell,\; \boldsymbol{k'}\in I, \boldk\neq\boldsymbol{k'}, \quad \bigcup_{\ell=1}^L I_\ell = I,
\end{equation}
are fulfilled.
The rank-1 lattice sizes $M_\ell$ are chosen distinctly from the set
\begin{align}
P^{I}_{\lambda,L_\mathrm{max}}&:=\left\{p_j\in P^I\colon p_j=\begin{cases}
\min\{p \in P^I\colon p>\lambda\}&\colon j=1\\
\min\{p \in P^I\colon p>p_{j-1}\}&\colon j=2,\ldots,L_\mathrm{max}.
\end{cases}
 \right\}\label{eq:def_PIlambdan}
\end{align}
of the $L_\mathrm{max}\in\N$ smallest prime numbers in
$$
P^I:=\{M'\in\N\colon M' \text{ prime with }|\{\zb k\bmod M'\colon \zb k\in I\}|=|I|\}
$$
larger than a certain $\lambda\in\N$.
With probability at least $1-\gamma$, $\Lambda$ returned by Algorithm~\ref{alg:construct_mr1l_I_distinct_primes:mod}
is a reconstructing multiple rank-1 lattice fulfilling properties~\eqref{equ:cond:reco_ml_mlfft2}.
For fixed oversampling factor $c>0$, the cardinality $M:=|\Lambda|$ of the multiple rank-1 lattice is in
\begin{align}
\mathcal{O}\big(\max\{|I|,N_I,-\log\gamma\,|\log|\log\gamma||\}(\log|I|-\log\gamma)\big),\label{eq:algo_mr1l_construction_sample_compl}
\end{align}
where $N_I:=\max_{j=1,\ldots, d}\{\max_{\zb k\in I}k_j-\min_{\zb l\in I}l_j\}$ is the expansion of $I$.
This follows since $M\lesssim L\max\{|I|, N_I, L\log L\}\lesssim \max\{|I|, N_I, L\log L\}(\log|I|-\log\gamma)$, cf. \cite[Corollary~3.7]{Kae17}, and
\begin{align*}
L\log L&\lesssim \log|I|\log\log|I|-\log\gamma\log|\log\gamma|-\log\gamma\log\log|I|+\log|I|\log|\log\gamma|\\
&\lesssim\log|I|\log\log|I|-\log\gamma\log|\log\gamma|\lesssim |I|-\log\gamma\,|\log|\log\gamma||.
\end{align*}
At this point, we stress on the fact that the oversampling factor $M/|I|$ does not depend on the dimension~$d$.
Moreover, the construction using Algorithm~\ref{alg:construct_mr1l_I_distinct_primes:mod} requires
\begin{align}
\mathcal{O}\big(\max\{|I|,N_I,-\log\gamma\,|\log|\log\gamma||\}(\log|I| + d + \log\log N_I)(\log|I|-\log\gamma)\big),\label{eq:algo_mr1l_construction_arithm_compl}
\end{align}
arithmetic operations.
This is caused by the complexity of determining the lattice sizes $M_\ell$, $\ell=1,\ldots,L_\mathrm{max}$, e.g.\ using the sieve of Eratosthenes, which is
\begin{align*}
&\lesssim M_{L_\mathrm{max}}\, \log\log M_{L_\mathrm{max}}\lesssim \max\{|I|, N_I, L\log L\}\,\log\log \max\{|I|, N_I, L\log L\}\\
& \lesssim \max\{|I|,N_I,-\log\gamma\,|\log|\log\gamma||\} (\log\log |I| + \log\log N_I + \log\log L)\\
& \lesssim \max\{|I|,N_I,-\log\gamma\,|\log|\log\gamma||\}\,(\log\log N_I)\, (\log\log |I| + \log\log L)\\
&\lesssim \max\{|I|,N_I,-\log\gamma\,|\log|\log\gamma||\}\,(\log\log N_I) \,(\log|I|-\log\gamma)
\end{align*}
and of the for loop, which is
\begin{align*}
L_\mathrm{max}(|I|\log|I|+d|I|)\lesssim |I|(\log|I|+d)(\log|I|-\log\gamma). 
\end{align*}

We remark that the proposed method represents a minor modification of \cite[Algorithm~4]{Kae17}.
After randomly choosing the generating vector $\boldz_L$ in line~\ref{alg:construct_mr1l_I_distinct_primes:mod:randomly_determine_z} of Algorithm~\ref{alg:construct_mr1l_I_distinct_primes:mod},
the original version in \cite[Algorithm~4]{Kae17} checks if the rank-1 lattice $\Lambda(\boldz_L,M_L)$ is able to reconstruct additional Fourier coefficients~$\hat{p}_\boldk$,
$\boldk\in I\setminus\tilde{I}$.
Otherwise, a new generating vector~$\boldz_L$ is randomly chosen and the check is repeated until at least one additional Fourier coefficient~$\hat{p}_\boldk$, $\boldk\in I\setminus\tilde{I}$,
can be reconstructed.
The analysis in \cite[Corollary~3.7]{Kae17} yields the same bounds
with respect to the success probability of finding a reconstructing multiple rank-1 lattice~$\Lambda$
and with respect to the cardinality~$|\Lambda|$
for \cite[Algorithm~4]{Kae17} and Algorithm~\ref{alg:construct_mr1l_I_distinct_primes:mod}.
However, the upper bound for the arithmetic complexity of \cite[Algorithm~4]{Kae17} is only with high probability, whereas the arithmetic complexity is guaranteed for Algorithm~\ref{alg:construct_mr1l_I_distinct_primes:mod} due to the realized modifications.
In practice, we observe that Algorithm~\ref{alg:construct_mr1l_I_distinct_primes:mod} may require slightly more samples than \cite[Algorithm~4]{Kae17}.

\begin{algorithm}[tb]
\caption{(Minor modification of \cite[Algorithm~4]{Kae17}). Determining reconstructing multiple rank-1 lattices with pairwise distinct lattice sizes,
fulfilling condition~\eqref{equ:cond:reco_ml_mlfft2} with high probability.
}\label{alg:construct_mr1l_I_distinct_primes:mod}
  \begin{tabular}{p{1.7cm}p{5.0cm}p{7.6cm}}
    Input: 	& $I\subset\Z^d$ 	& frequency index set\\
    		& $c\in(1,\infty)\subset\R$ 	& oversampling factor\\
    		& $\gamma\in(0,1)\subset\R$			& upper bound on failure probability
  \end{tabular}
		
  \begin{algorithmic}[1]
  	\STATE Determine maximal number of rank-1 lattices $L_\mathrm{max}:=\ceil{c^2/(c-1)^2 \, (\ln |I|-\ln\gamma)/2}$.
  	\STATE Determine lower bound $\lambda:=c\;(|I|-1)$ for lattice sizes.
	\STATE Determine set $P^{I}_{\lambda,L_\mathrm{max}}$ of prime lattice sizes, cf. \eqref{eq:def_PIlambdan}, and arrange $p_1<\ldots<p_{L_\mathrm{max}}$.
	\STATE Initialize index set $\tilde{I}:=\emptyset$ of frequencies belonging to reconstructible Fourier coefficients.
	\STATE {\bfseries for} {$L:=1,\ldots,L_\mathrm{max}$} {\bfseries do}
		\STATE \hspace{1em} Choose lattice size $M_{L}=p_{L}\in P^{I}_{\lambda,L_\mathrm{max}}$.\label{alg:construct_mr1l_I_distinct_primes:determinstic_not_random_lattice_sizes}
		\STATE \hspace{1em} Choose generating vector $\zb z_{L}$ from $[0,M_{L}-1]^d\cap\Z^d$ uniformly at random.\label{alg:construct_mr1l_I_distinct_primes:mod:randomly_determine_z}
	    \STATE \hspace{1em} Determine $I_L:=\{\zb k\in I\colon \not\exists \zb h\in I\setminus\{\zb k\}\text{ with }\zb k\cdot \zb z_{L}\equiv\zb h\cdot \zb z_{L}\imod{M_{L}}\}$, cf.~\eqref{equ:cond:reco_ml_mlfft2}.
	    \STATE \hspace{1em} Update $\tilde{I}:=\tilde{I}\cup I_L$.
	    \STATE \hspace{1em} If $|\tilde{I}|=|I|$, then exit {\bfseries for} loop.
	\STATE {\bfseries end for} $L$
  \end{algorithmic}
  \begin{tabular}{p{1.7cm}p{5.0cm}p{7.6cm}}
    Output: & $M_1,\ldots,M_{L}$ & lattice sizes of rank\mbox{-}1 lattices and\\
	    & $\zb z_1,\ldots,\zb z_{L}$ & generating vectors of rank\mbox{-}1 lattices such that\\
	    & $\Lambda(\zb z_1,M_1,\ldots,\zb z_{L},M_{L})$ &  is a reconstructing multiple rank\mbox{-}1 lattice for~$I$ 
	   	 with probability at least $1-\gamma$  
\\    \cmidrule{1-3}
 \end{tabular}
  \begin{tabular}{p{1.7cm}p{5.0cm}p{7.6cm}}
    Complexity: & $\mathcal{O}\big(|I|\,(\log|I|+d)\log|I|\big)$&{for $|I|\gtrsim N_I$, fixed $c$ and fixed $\gamma$, where \newline $N_I:=\max_{j=1,\ldots, d}\{\max_{\zb k\in I}k_j-\min_{\zb l\in I}l_j\}$ \newline is the expansion of $I$.} 
  \end{tabular}
\end{algorithm}

Another construction method is \cite[Algorithm~6]{Kae17}, which searches for a reconstructing multiple rank-1 lattice such that the properties
$$
 \boldk\cdot\boldz_\ell\not\equiv \boldsymbol{k'}\cdot\boldz_\ell\imod{M_\ell} \text{ for all } \boldk\in I_\ell,\; \boldsymbol{k'}\in I\setminus\left\{\bigcup_{\ell'=1}^{\ell-1} I_{\ell'}\right\},\;\boldk\neq\boldsymbol{k'},\quad \bigcup_{\ell=1}^L I_\ell = I,
$$
are fulfilled. Under mild assumptions, \cite[Algorithm~6]{Kae17} returns a reconstructing multiple rank-1 lattice
of cardinality $\mathcal{O}(|I|\log^2|I|)$ with high probability, where the constants do not depend on the dimension $d$, and requires $\mathcal{O}(|I| (d+\log|I|) \log^3|I|)$ arithmetic operations with high probability.

\begin{algorithm}[tp]
\caption{Non-iterative reconstruction of a trigonometric polynomial $p$ from sampling values along reconstructing multiple rank-1 lattices for known frequency index set $I$.}\label{algo:mlfft2:ifft_direct}
  \begin{tabular}{p{1.7cm}p{5.0cm}p{7.6cm}}
    Input:      & $I\subset\Z^d$ & frequency index set, $|I|<\infty$ \\
                & $\Lambda:=\Lambda(\boldz_1,M_1,\ldots,\boldz_L,M_L)$ & reconstructing multiple rank-1 lattice for $I$ fulfilling condition~\eqref{equ:cond:reco_ml_mlfft2} \\
                & $\left(p(\boldsymbol{\tilde{x}}_j)\right)_{\boldsymbol{\tilde{x}}_j\in\Lambda}$  & sampling values of trigonometric polynomial~$p$ \\
  \end{tabular}
  \begin{algorithmic}
     \STATE Initialize $\mathtt{counter}[\boldk]:=0$ and $\tilde{\hat{p}}_\boldk:=0$ for $\boldk\in I$.
     \STATE {\bfseries for} $\ell:=1,\ldots,L$ {\bfseries do}
     \STATE \hspace{1em} Determine $I_\ell:=\{\zb k\in I\colon \not\exists \zb h\in I\setminus\{\zb k\}\text{ with }\zb k\cdot \zb z_\ell\equiv\zb h\cdot \zb z_\ell\imod{M_\ell}\}$, cf.~\eqref{equ:cond:reco_ml_mlfft2}.
     \STATE \hspace{1em} Compute $\tilde{\hat{g}}_\boldk:=\frac{1}{M_\ell} \sum_{j=0}^{M_\ell-1} p\left(\frac{j}{M_\ell}\boldz_\ell \bmod \boldone\right) \,\mathrm{e}^{-2\pi\mathrm{i} j\boldk\cdot\boldz_\ell /M_\ell}$ for $\boldk\in I_\ell$, using inverse rank-1 lattice FFT, cf.\ \cite[Algorithm 3.2]{kaemmererdiss}.
     \STATE \hspace{1em} Set $\mathtt{counter}[\boldk]:=\mathtt{counter}[\boldk]+1$ and $\tilde{\hat{p}}_\boldk:=\tilde{\hat{p}}_\boldk+\tilde{\hat{g}}_\boldk$ for $\boldk\in I_\ell$.
     \STATE {\bfseries end for} $\ell$
     \STATE Set $\tilde{\hat{p}}_\boldk:=\tilde{\hat{p}}_\boldk/\mathtt{counter}[\boldk]$ for $\boldk\in I$.
  \end{algorithmic}
  \begin{tabular}{p{1.7cm}p{5.1cm}p{7.5cm}}
    Output: & $\left(\tilde{\hat{p}}_\boldk\right)_{\boldk\in I}$ & reconstructed Fourier coefficients \\
    \cmidrule{1-3}
 \end{tabular}
  \begin{tabular}{p{1.75cm}p{5.05cm}p{7.5cm}}
    Complexity: & $\mathcal{O}\big(\widetilde{M}\log \widetilde{M} + L\,|I|\,(d + \log |I|)\big)$ & $|\Lambda|\leq \sum_{\ell=1}^{L}M_\ell =: \widetilde{M}$\\
  \end{tabular}
\end{algorithm}

Besides the small cardinalities and fast construction algorithms, a further main advantage of reconstructing multiple rank-1 lattices is the existence of a direct and fast inversion method for computing Fourier coefficients $\hat{p}_\boldk$ from sampling values,
cf. Algorithm~\ref{algo:mlfft2:ifft_direct} for multiple rank-1 lattices constructed by \cite[Algorithm~4]{Kae17} or Algorithm~\ref{alg:construct_mr1l_I_distinct_primes:mod} as well as cf. \cite[Algorithm~6]{Kae16} for multiple rank-1 lattices constructed by \cite[Algorithm~6]{Kae17}.

\subsection{Dimension-incremental method}
\label{sec:pre:reco_incr}

We use the dimension-incremental approach from~\cite{PoVo14}. This approach proceeds similarly as a dimension-incremental method for anharmonic trigonometric polynomials based on Prony's method in \cite{PoTa12}.
For a given search domain $\Gamma\subset\Z^d$, $|\Gamma|<\infty$, which may be very large like e.g.\ a full grid~$\hat{G}_N^d$, a frequency index set $I\subset\Gamma$ containing the approximately largest Fourier coefficients of a function~$f$ under consideration is determined based on samples. In doing so, the index set $I$ is constructed component-by-component in a dimension-incremental way starting with the first component.

In order to describe the method, we introduce additional notation from~\cite{PoVo14} and we assume that the function under consideration is a multivariate trigonometric polynomial~$p$.
We denote the projection of a frequency $\boldk:=(k_1,\ldots,k_d)^\top\in\Z^d$
to the components $\boldi:=(i_1,\ldots,i_m)\in\left\{\boldsymbol{\iota} \in \{1,\ldots,d\}^m\colon \iota_t\neq\iota_{t'} \text{ for }t\neq t' \right\}$
by $\mathcal{P}_\boldi(\boldk) := (k_{i_1},\ldots,k_{i_m})^\top\in\Z^m$.
Correspondingly, we define the projection of a frequency index set
$I\subset\Z^d$ to the components $\boldi$ by
$\mathcal{P}_\boldi(I):=\{(k_{i_1},\ldots,k_{i_m})\colon \boldk\in I\}$.
Using this notation, the general approach is the following:
\begin{enumerate}
 \item \label{enum:reco:incr_general:I1}
         Determine the first components of the unknown frequency index set, i.e., determine an index set $I^{(1)}\subseteq\mathcal{P}_1(\Gamma)$
         which should be identical to the projection $\mathcal{P}_1(\mathrm{supp}\,\hat{p})$ or contain this projection,
         $I^{(1)}\supseteq\mathcal{P}_1(\mathrm{supp}\,\hat{p})$.
 \item \label{enum:reco:incr_general:I2}
       For dimension increment step $t=2,\ldots,d$, i.e., for each additional dimension:
 \begin{enumerate}
   \item \label{enum:reco:incr_general:I2:I_t}
         Determine the $t$-th components of the unknown frequency index set, i.e., determine an index set $I^{(t)}\subseteq\mathcal{P}_t(\Gamma)$ 
         which should be identical to the projection $\mathcal{P}_t(\mathrm{supp}\,\hat{p})$ or contain this projection,
         $I^{(t)}\supseteq\mathcal{P}_t(\mathrm{supp}\,\hat{p})$.
   \item \label{enum:reco:incr_general:I2:X_1_t}
         Determine a suitable sampling set $\mathcal{X}^{(1,\ldots,t)}\subset\T^d$, $\vert\mathcal{X}^{(1,\ldots,t)}\vert\ll\vert\Gamma\vert$, which allows to detect those frequencies from the index set $(I^{(1,\ldots,t-1)} \times I^{(t)})\cap\mathcal{P}_{(1,\ldots,t)}(\Gamma)$ belonging to non-zero Fourier coefficients $\hat{p}_\boldk$.
   \item \label{enum:reco:incr_general:I2:sample}
         Sample the trigonometric polynomial $p$ along the nodes of the sampling set $\mathcal{X}^{(1,\ldots,t)}$.
   \item \label{enum:reco:incr_general:I2:p_hat_tilde} Compute the Fourier coefficients $\tilde{\hat{p}}_{(1,\ldots,t),\boldk}$, $\boldk\in (I^{(1,\ldots,t-1)} \times I^{(t)})\cap\mathcal{P}_{(1,\ldots,t)}(\Gamma)$.
   \item \label{enum:reco:incr_general:I2:non-zero}
       Determine the non-zero Fourier coefficients from $\tilde{\hat{p}}_{(1,\ldots,t),\boldk}$, $\boldk\in (I^{(1,\ldots,t-1)} \times I^{(t)})\cap\mathcal{P}_{(1,\ldots,t)}(\Gamma)$,
       and obtain the index set $I^{(1,\ldots,t)}$ of detected frequencies.
       The $I^{(1,\ldots,t)}$ index set should be equal to the projection $\mathcal{P}_{(1,\ldots,t)}(\mathrm{supp}\,\hat{p})$.
 \end{enumerate}
 \item \label{enum:reco:incr_general:I3}
       Use the index set $I^{(1,\ldots,d)}$ and the computed Fourier coefficients $\tilde{\hat{p}}_{(1,\ldots,d),\boldk}$, $\boldk\in I^{(1,\ldots,d)}$,
       as an approximation for the support $\mathrm{supp}\,\hat{p}$ and the Fourier coefficients $\hat{p}_\boldk$, $\boldk\in\mathrm{supp}\,\hat{p}$.
\end{enumerate}

The proposed approach includes the construction of a suitable sampling set in step~\ref{enum:reco:incr_general:I2:X_1_t} and the computation of (projected) Fourier coefficients in step~\ref{enum:reco:incr_general:I2:p_hat_tilde}.
There exist different methods for the realization of these steps.

In~\cite{PoVo14}, the utilization of reconstructing rank-1 lattices was mainly considered as sampling sets, which allows for efficiently computing the (projected) Fourier coefficients using a so-called rank-1 lattice FFT, i.e., a single one-dimensional FFT followed by a simple index transform, cf.\ \cite[Algorithm 3.2]{kaemmererdiss}. This approach yielded high stability and reliability in numerical tests, but required a relatively high number of samples and the construction of reconstructing rank-1 lattices may involve high computational costs.

Alternatively, for computing the (projected) Fourier coefficients, sub-sampling along reconstructing rank-1 lattices using $\ell_1$ minimization or using a variant of Prony's method as well as
randomly obtained generated sets using $\ell_1$ minimization were also considered in~\cite{PoVo14}. These approaches performed well when determining the unknown frequency index sets~$I$ and Fourier coefficients~$\hat{p}_\boldk$ of sparse multivariate trigonometric polynomials $p$ based on unperturbed sampling values. However, for the case of noisy sampling values or functions $f\in L_1(\T^d)\cap\mathcal{C}(\T^d)$ with infinite support in frequency domain, these methods failed in numerical tests.

In the next section, we present a novel approach, which is based on using multiple rank-1 lattices, which does not have these limitations.

\section{Method}
\label{sec:method}

In this section, we 
present a substantial improvement of the dimension-incremental approach developed in \cite{PoVo14}.
On the one hand, we change the sampling sets to multiple rank-1 lattices $\Lambda=\Lambda(\boldz_1,M_1,\ldots,\boldz_L,M_L)$,
which significantly reduces the number of required sampling values compared to single rank-1 lattices~$\Lambda(\boldz,M)$.
On the other hand, we describe modifications on \cite[Algorithm~2]{PoVo14}, that allow for estimates of the success probability,
the total number of required sampling values, and the overall complexity of the resulting algorithms.
In particular for sparse trigonometric polynomials, these modifications will entail 
detailed estimates that only depend on properties of the frequency support of the sparse trigonometric polynomial
and a factor that is logarithmic in the failure probability, cf. Section~\ref{sec:theoretical}.
  
In Section~\ref{sec:method:algorithm}, we present two algorithms based on sampling along multiple rank-1 lattices. These algorithms perform the steps as generally explained in Section~\ref{sec:pre:reco_incr}.
Afterwards, in Section~\ref{sec:method:complexity}, we discuss the number of required samples and arithmetic operations of the proposed methods in general.

\subsection{Algorithm}
\label{sec:method:algorithm}

We present the first proposed realization of the dimension-incremental approach from Section~\ref{sec:pre:reco_incr} as Algorithm~\ref{algo:sfft_mlfft},
which uses multiple rank-1 lattices as sampling nodes, cf. Section~\ref{sec:pre:multiple_r1l}.

Compared to \cite[Algorithm~2]{PoVo14}, we perform the following modifications. In step~\ref{enum:reco:incr_general:I2:X_1_t}, we apply Algorithm~\ref{alg:construct_mr1l_I_distinct_primes:mod}
on the frequency index set $(I^{(1,\ldots,t-1)} \times I^{(t)})\cap\mathcal{P}_{(1,\ldots,t)}(\Gamma)$ in order to build a reconstructing multiple \hbox{rank-1} lattice. This change 
of the construction of the sampling sets drastically reduces the arithmetic complexity of the overall algorithm as discussed in Section~\ref{sec:method:complexity}.
Additionally, the resulting sampling set $\mathcal{X}^{(1,\ldots,t)}$ consists of distinctly less nodes and, consequently, the number of sampling values in step~\ref{enum:reco:incr_general:I2:sample} is reduced.
Afterwards in step~\ref{enum:reco:incr_general:I2:p_hat_tilde}, the projected Fourier coefficients $\tilde{\hat{p}}_{(1,\ldots,t),\boldk}$, $\boldk\in (I^{(1,\ldots,t-1)} \times I^{(t)})\cap\mathcal{P}_{(1,\ldots,t)}(\Gamma)$,
are computed using Algorithm~\ref{algo:mlfft2:ifft_direct}. 
Moreover, we do not need to perform the additional step~2f from \cite[Algorithm~2]{PoVo14}, which also reduces the number of arithmetic operations.

Another essential modification compared to the algorithms presented in \cite{PoVo14}
concerns the thresholding behavior.
In \cite{PoVo14} a relative threshold parameter $\theta\in\R^+$ with respect to a currently numerically computed Fourier coefficient $\tilde{\hat{p}}_{(1,\ldots,t),\boldk}$ was used.
Here, we utilize an absolute threshold $\delta\in\R^+$, which is an input parameter of Algorithm~\ref{algo:sfft_mlfft}. This absolute threshold is required for our theoretical considerations in Section~\ref{sec:theoretical}, where we show estimates for failure probabilities and for the number of required detection iterations $r\in\N$.

Moreover, we introduce an additional input parameter $b\in\N$, which is the maximal number of multiple rank-1 lattice searches per component for the search of reconstructing multiple rank-1 lattices.
This parameter $b$ is required in order to achieve a deterministic runtime behavior and
is considered theoretically in Lemma~\ref{lem:reco_mr1l:rep_b} of Section~\ref{sec:theoretical:prob}.
In practice, choosing the parameters $r$ and $b$ as a small constant  (e.g. $\leq 10$) should often suffice as observed in the numerical tests in Section~\ref{sec:numerics}.

\begin{algorithm}[htp]
\caption{Reconstruction of a multivariate trigonometric polynomial $p$ from sampling values along multiple rank-1 lattices.}\label{algo:sfft_mlfft}
  \begin{tabular}{p{1.2cm}p{2.0cm}p{11.1cm}}
    Input:  
                & $\Gamma\subset\Z^d$ \hfill & search space in frequency domain, superset for $\mathrm{supp}\,\hat{p}$ \\
                & $p(\circ)$ & trigonometric polynomial $p$ as black box (function handle) \\
                & $\delta\in\R^+$ & absolute threshold \\ 
                & $s,s_\mathrm{local}\in\N$ & sparsity parameter ($s_\mathrm{local}:=s$ by default) \\ 
                & $r\in\N$ & number of detection iterations (e.g. see Theorem~\ref{thm:r_complexity} or Section~\ref{sec:numerics}) \\
                & $b\in\N$ & maximal number of multiple rank-1 lattice searches per component
  \end{tabular}
  \vspace{-0.6em}
  \begin{algorithmic}
   \item[(step \ref{enum:reco:incr_general:I1})]
   \STATE Set $K_1:=\max(\mathcal{P}_1(\Gamma))-\min(\mathcal{P}_1(\Gamma))+1$, $I^{(1)}:=\emptyset$.
   \STATE {\bfseries for} $i:=1,\ldots,r$ {\bfseries do}
   \STATE \hspace{1em} Choose components $x_2',\ldots,x_d'\in\T$ of sampling nodes uniformly at random.
   \STATE \hspace{1em} Compute $\tilde{\hat{p}}_{1,k_1}:=\frac{1}{K_1} \sum_{\ell=0}^{K_1-1} p\left(\frac{\ell}{K_1},x_2',\ldots,x_d'\right) \,\mathrm{e}^{-2\pi\mathrm{i} \ell k_1 / K_1}$, $k_1\in\mathcal{P}_1(\Gamma)$, with FFT.
   \STATE \hspace{1em} Set $I^{(1)}:=I^{(1)}\cup\{k_1\in\mathcal{P}_1(\Gamma)\colon \text{(up to) $s_\mathrm{local}$-largest values } \vert\tilde{\hat{p}}_{1,k_1}\vert \geq \delta 
   \}$.
   \STATE {\bfseries end for} $i$
   \item[] \vspace{-0.9em}
   \item[(step \ref{enum:reco:incr_general:I2})]
   \item[{\bfseries for} $t:=2,\ldots,d$ {\bfseries do}]
     \item[(step \ref{enum:reco:incr_general:I2:I_t})]
     \STATE Set $K_t:=\max(\mathcal{P}_t(\Gamma))-\min(\mathcal{P}_t(\Gamma))+1$, $I^{(t)}:=\emptyset$.
     \STATE {\bfseries for} $i:=1,\ldots,r$ {\bfseries do}
     \STATE \hspace{1em} Choose components $x_1',\ldots,x_{t-1}',x_{t+1}',\ldots,x_d'\in\T$ of sampling nodes uniformly at random.
     \STATE \hspace{1em} $\tilde{\hat{p}}_{t,k_t}:=
     \sum\limits_{\ell=0}^{K_t-1} p\left(x_1',\ldots,x_{t-1}',\frac{\ell}{K_t},x_{t+1}',\ldots,x_d'\right) \mathrm{e}^{-2\pi\mathrm{i}\ell k_t/K_t}$, $k_t\in\mathcal{P}_t(\Gamma)$, using FFT.
     \STATE \hspace{1em} Set $I^{(t)}:=I^{(t)}\cup\{k_t\in\mathcal{P}_t(\Gamma)\colon \text{(up to) $s_\mathrm{local}$-largest values } \vert\tilde{\hat{p}}_{t,k_t}\vert \geq \delta 
     \}$.
     \STATE {\bfseries end for} $i$
     \item[(step \ref{enum:reco:incr_general:I2:X_1_t})]
     \STATE If $t<d$, set $\tilde r:=r$ and $\tilde s:=s_\mathrm{local}$, otherwise $\tilde r:=1$ and $\tilde s:=s$. 
     Set $I^{(1,\ldots,t)}:=\emptyset$.
     \STATE Search for reconstructing multiple rank-1 lattice $\Lambda(\boldz_1,M_1,\ldots,\boldz_L,M_L)$ for $(I^{(1,\ldots,t-1)} \times I^{(t)})\cap\mathcal{P}_{(1,\ldots,t)}(\Gamma)$ using {\bfseries Algorithm~\ref{alg:construct_mr1l_I_distinct_primes:mod}} with algorithm parameters $c:=2$ and $\gamma:=0.5$; repeat up to $b-1$ times or until reconstruction property~\eqref{equ:cond:reco_ml_mlfft2} is fulfilled.
     \newline $\Lambda:=\Lambda(\boldz_1,M_1,\ldots,\boldz_L,M_L)$ consists of nodes $\boldsymbol{\tilde{x}}_j\in\T^t$, $j=0,\ldots,M^{(t)}-1$, $M^{(t)}=|\Lambda|$. 
     \STATE {\bfseries for} $i:=1,\ldots,\tilde r$ {\bfseries do}
     \STATE \hspace{1em} Choose components $x_{t+1}',\ldots,x_d'\in\T$ of sampling nodes uniformly at random.
     \STATE \hspace{1em} Set $\mathcal{X}^{(1,\ldots,t)}:=\{\boldx_j:=(\boldsymbol{\tilde{x}}_j,x_{t+1}',\ldots,x_d') \bmod \boldone\colon j=0,\ldots,M^{(t)}-1\}\subset\T^d$.
     \item[(step \ref{enum:reco:incr_general:I2:sample})] Sample $p$ along the nodes of the sampling set $\mathcal{X}^{(1,\ldots,t)}$.
     \item[(step \ref{enum:reco:incr_general:I2:p_hat_tilde})]
     Compute $\tilde{\hat{p}}_{(1,\ldots,t),\boldk}$
     for $\boldk\in (I^{(1,\ldots,t-1)} \times I^{(t)})\cap\mathcal{P}_{(1,\ldots,t)}(\Gamma)$ from sampling values~$p(\boldx_j)$, $\boldx_j\in\mathcal{X}^{(1,\ldots,t)}$, with inverse multiple rank-1 lattice FFT
     using {\bfseries Algorithm~\ref{algo:mlfft2:ifft_direct}}.
     \\
     \item[(step \ref{enum:reco:incr_general:I2:non-zero})]
     Set $I^{(1,\ldots,t)}:=I^{(1,\ldots,t)}\cup \{ \boldk\in (I^{(1,\ldots,t-1)} \times I^{(t)})\cap\mathcal{P}_{(1,\ldots,t)}(\Gamma)\colon \newline \text{(up to) $\tilde s$-largest values } \vert\tilde{\hat{p}}_{(1,\ldots,t),\boldk}\vert \geq \delta 
     \}$.
     \STATE {\bfseries end for} $i$
   \item[{\bfseries end for} $t$]
  \end{algorithmic}
\end{algorithm}

\begin{algorithm*}
  \ContinuedFloat
  \caption{continued.}
  \begin{algorithmic}
   \item[(step \ref{enum:reco:incr_general:I3})] Set $I:=I^{(1,\ldots,d)}$ and $\boldsymbol{\tilde{\hat{p}}}:=\left(\tilde{\hat{p}}_{(1,\ldots,d),\boldk}\right)_{\boldk\in I}$.
  \end{algorithmic}
  \begin{tabular}{p{1.2cm}p{2.0cm}p{11.1cm}}
    Output: & $I\subset\Gamma\subset\Z^d$ & index set of detected frequencies, $|I|\leq \min \{s,|\Gamma|\}$ \\
            & $\boldsymbol{\tilde{\hat{p}}}\in\C^{\vert I\vert}$ & corresponding Fourier coefficients, $|\tilde{\hat{p}}_{(1,\ldots,d),\boldk}| \geq \delta$ \\
  \end{tabular}
\end{algorithm*}

Additionally,
we present a second dimensional-incremental reconstruction algorithm.
To this end, we modify Algorithm~\ref{algo:sfft_mlfft} by using \cite[Algorithm~6]{Kae17} for the determination of the reconstructing multiple rank-1 lattices and obtain Algorithm~\ref{algo:sfft_mlfft6}.
We stress on the fact that \cite[Algorithm~6]{Kae17} does not ensure to terminate. However,
if it is terminating, the output is always a reconstructing multiple rank-1 lattice.
Therefore, we omit a repeatedly invocation of \cite[Algorithm~6]{Kae17} in step 2b and we remove the parameter $b$.
Since the construction approach of the multiple rank-1 lattices changes, we need to apply
the corresponding FFT algorithms, i.e.,
the projected Fourier coefficients $\tilde{\hat{p}}_{(1,\ldots,t),\boldk}$, $\boldk\in (I^{(1,\ldots,t-1)} \times I^{(t)})\cap\mathcal{P}_{(1,\ldots,t)}(\Gamma)$,
in step~\ref{enum:reco:incr_general:I2:p_hat_tilde} are now computed using the inverse multiple rank-1 lattice FFT from \cite[Algorithm~6]{Kae16}.
As discussed in Section~\ref{sec:method:complexity}, the sample and arithmetic complexities of Algorithm~\ref{algo:sfft_mlfft6} are slightly higher by a logarithmic factor and non-deterministic.
However, in the numerical tests in Section~\ref{sec:numerics}, the observed total number of samples and runtimes are usually smaller for Algorithm~\ref{algo:sfft_mlfft6}.

Moreover, we will apply parts of the theoretical results from Section~\ref{sec:theoretical} to
the dimension-incremental approach using single rank-1 lattices as sampling sets.
For this, we need to modify the thresholding behavior of \cite[Algorithm~2]{PoVo14}.
As described above, we use an absolute threshold $\delta\in\R^+$ as input parameter
and this yields Algorithm~\ref{algo:sfft_a2r1labs}.

\begin{algorithm}[t]
\caption{(Modified Algorithm~\ref{algo:sfft_mlfft}). Reconstruction of a multivariate trigonometric polynomial $p$ from sampling values along reconstructing multiple rank-1 lattices.}\label{algo:sfft_mlfft6}
  \begin{tabular}{p{1.2cm}p{2.0cm}p{11.1cm}}
    Input:  
                & $\Gamma\subset\Z^d$ \hfill & search space in frequency domain, superset for $\mathrm{supp}\,\hat{p}$ \\
                & $p(\circ)$ & trigonometric polynomial $p$ as black box (function handle) \\
                & $\delta\in\R^+$ & absolute threshold \\ 
                & $s,s_\mathrm{local}\in\N$ & sparsity parameter ($s_\mathrm{local}:=s$ by default) \\ 
                & $r\in\N$ & number of detection iterations \\
  \end{tabular}
  \vspace{-0.6em}
  \begin{algorithmic}
   \item[(step \ref{enum:reco:incr_general:I1})]
   \item[\vdots]
   \item[(step \ref{enum:reco:incr_general:I2})]
   \item[{\bfseries for} $t:=2,\ldots,d$ {\bfseries do}]
   \item[\hspace{0.5cm}\vdots]
     \item[(step \ref{enum:reco:incr_general:I2:X_1_t})]
     \STATE If $t<d$, set $\tilde r:=r$ and $\tilde s:=s_\mathrm{local}$, otherwise $\tilde r:=1$ and $\tilde s:=s$. 
     Set $I^{(1,\ldots,t)}:=\emptyset$.
     \STATE Search for reconstructing multiple rank-1 lattice $\Lambda(\boldz_1,M_1,\ldots,\boldz_L,M_L)$ for $(I^{(1,\ldots,t-1)} \times I^{(t)})\cap\mathcal{P}_{(1,\ldots,t)}(\Gamma)$ using {\bfseries\cite[Algorithm~6]{Kae17}} with algorithm parameters $c:=2$ and $\delta':=0.5$ 
    \STATE {\bfseries for} $i:=1,\ldots,\tilde r$ {\bfseries do}
   \item[\hspace{1cm}\vdots]
      \item[(step \ref{enum:reco:incr_general:I2:p_hat_tilde})]
     Compute $\tilde{\hat{p}}_{(1,\ldots,t),\boldk}$
     for $\boldk\in (I^{(1,\ldots,t-1)} \times I^{(t)})\cap\mathcal{P}_{(1,\ldots,t)}(\Gamma)$ from sampling values~$p(\boldx_j)$, $\boldx_j\in\mathcal{X}^{(1,\ldots,t)}$, with inverse multiple rank-1 lattice FFT
     using {\bfseries\cite[Algorithm~6]{Kae16}}.
     \\
     \item[\hspace{1cm}\vdots]
   \STATE {\bfseries end for} $i$
   \item[{\bfseries end for} $t$]
   \item[(step \ref{enum:reco:incr_general:I3})]
   \STATE Compute $\tilde{\hat{p}}_{(1,\ldots,d),\boldk}$
     for $\boldk\in I^{(1,\ldots,d)}$ from sampling values~$p(\boldx_j)$, $\boldx_j\in\mathcal{X}^{(1,\ldots,d)}$, with inverse multiple rank-1 lattice FFT
     using \cite[Algorithm~6]{Kae16}.
   \STATE Set $I:=I^{(1,\ldots,d)}$ and $\boldsymbol{\tilde{\hat{p}}}:=\left(\tilde{\hat{p}}_{(1,\ldots,d),\boldk}\right)_{\boldk\in I}$.
  \end{algorithmic}
  \begin{tabular}{p{1.2cm}p{2.0cm}p{11.1cm}}
    Output: & $I\subset\Gamma\subset\Z^d$ & index set of detected frequencies, $|I|\leq \min \{s,|\Gamma|\}$ \\
            & $\boldsymbol{\tilde{\hat{p}}}\in\C^{\vert I\vert}$ & corresponding Fourier coefficients, $|\tilde{\hat{p}}_{(1,\ldots,d),\boldk}| \geq \delta$ \\
  \end{tabular}
\end{algorithm}

\begin{algorithm}[htp]
\caption{(Modified version of \cite[Algorithm~2]{PoVo14}). Reconstruction of a multivariate trigonometric polynomial $p$ from sampling values along (single) reconstructing rank-1 lattices.}\label{algo:sfft_a2r1labs}
  \begin{tabular}{p{1.2cm}p{2.0cm}p{11.1cm}}
    Input:  
                & $\Gamma\subset\Z^d$ \hfill & search space in frequency domain, superset for $\mathrm{supp}\,\hat{p}$ \\
                & $p(\circ)$ & trigonometric polynomial $p$ as black box (function handle) \\
                & $\delta\in\R^+$ & absolute threshold \\ 
                & $s,s_\mathrm{local}\in\N$ & sparsity parameter ($s_\mathrm{local}:=s$ by default) \\ 
                & $r\in\N$ & number of detection iterations (e.g. see Theorem~\ref{cor:r_complexity:single_r1l} or Section~\ref{sec:numerics})
  \end{tabular}
  \vspace{-0.6em}
  \begin{algorithmic}
   \item[(step \ref{enum:reco:incr_general:I1})]
   \item[\vdots]
   \STATE \textbf{for} $i:=1,\ldots,r$ \textbf{do}
     \item[\hspace{1cm}\vdots]
     \STATE \hspace{1em} $I^{(1)}:=I^{(1)}\cup\{k_1\in\mathcal{P}_1(\Gamma)\colon \text{(up to) $s_\mathrm{local}$-largest values } \vert\tilde{\hat{p}}_{1,k_1}\vert \geq \delta \}$
   \STATE \textbf{end}~\textbf{for} $i$
   \item[\vdots]
   \item[(step \ref{enum:reco:incr_general:I2})]
   \item[\textbf{for} $t:=2,\ldots,d$ \textbf{do}]
     \item[\hspace{0.5cm}\vdots]
     \STATE \textbf{for} $i:=1,\ldots,r$ \textbf{do}
      \item[\hspace{1cm}\vdots]
     \STATE \hspace{1em} Set $I^{(t)}:=I^{(t)}\cup\{k_t\in\mathcal{P}_t(\Gamma)\colon \text{(up to) $s_\mathrm{local}$-largest values } \vert\tilde{\hat{p}}_{t,k_t}\vert \geq \delta \}$.
     \STATE \textbf{end}~\textbf{for} $i$
     \item[(step \ref{enum:reco:incr_general:I2:X_1_t})]
     \STATE If $t<d$, set $\tilde r:=r$ and $\tilde s:=s_\mathrm{local}$, otherwise $\tilde r:=1$ and $\tilde s:=s$. 
     Set $I^{(1,\ldots,t)}:=\emptyset$.
     \item[\hspace{0.5cm}\vdots]
     \STATE \textbf{for} $i:=1,\ldots,\tilde r$ \textbf{do}
      \item[\hspace{1cm}\vdots]
     \item[(step \ref{enum:reco:incr_general:I2:non-zero})]
     \STATE \hspace{1em} Set $I^{(1,\ldots,t)}:=I^{(1,\ldots,t)}\cup \{ \boldk\in (I^{(1,\ldots,t-1)} \times I^{(t)})\cap\mathcal{P}_{(1,\ldots,t)}(\Gamma)\colon \newline \text{(up to) $\tilde s$-largest values } \vert\tilde{\hat{p}}_{(1,\ldots,t),\boldk}\vert \geq  \delta \}$.
     \STATE \textbf{end}~\textbf{for} $i$
   \item[\hspace{0.5cm}\vdots]
   \item[\textbf{end} \textbf{for} $t$]
   \item[(step \ref{enum:reco:incr_general:I3})]
   Set $I:=I^{(1,\ldots,d)}$ and $\boldsymbol{\tilde{\hat{p}}}:=\left(\tilde{\hat{p}}_{(1,\ldots,d),\boldk}\right)_{\boldk\in I}$.
  \end{algorithmic}
  \begin{tabular}{p{1.2cm}p{2.0cm}p{11.1cm}}
    Output: & $I\subset\Gamma\subset\Z^d$ & index set of detected frequencies, $|I|\leq \min \{s,|\Gamma|\}$ \\
            & $\boldsymbol{\tilde{\hat{p}}}\in\C^{\vert I\vert}$ & corresponding Fourier coefficients, $|\tilde{\hat{p}}_{(1,\ldots,d),\boldk}| \geq \delta$ \\
  \end{tabular}
\end{algorithm}

\begin{Remark}
The dimension-incremental method as described in Section~\ref{sec:pre:reco_incr} and realized by \cite[Algorithm~1,2]{PoVo14} as well as Algorithm~\ref{algo:sfft_mlfft}, \ref{algo:sfft_mlfft6} and~\ref{algo:sfft_a2r1labs} proceeds component-by-component starting with the first dimension, then the second et cetera. In principle, any order / permutation of the components $\{1,2,\ldots,d\}$ may be used. Depending on the actual function under consideration, a different order may be beneficial or not, since this may influence the aliasing behavior in the various dimension increment steps $t$.
\end{Remark}

\begin{Remark}\label{rem:deterministic}
As discussed in \cite[Section~2.2.1]{PoVo14}, we may set the number of detection iterations $r:=1$ for \cite[Algorithm~1,2]{PoVo14}, Algorithm~\ref{algo:sfft_mlfft}, Algorithm~\ref{algo:sfft_mlfft6} and Algorithm~\ref{algo:sfft_a2r1labs}
if the Fourier coefficients $\hat{p}_\boldk$ of the trigonometric polynomial $p$ under consideration fulfill a certain property,
namely, that the signs of the real part $\mathrm{Re}(\hat{p}_\boldk)$ of all non-zero Fourier coefficients $\hat{p}_\boldk$ are identical $\neq 0$ or that the signs of the imaginary part $\mathrm{Im}(\hat{p}_\boldk)$ are identical $\neq 0$. Then, we do not choose the components $x_1',\ldots,x_d'\in\T$ randomly but fix them to zero. In the case of \cite[Algorithm~1,2]{PoVo14} and Algorithm~\ref{algo:sfft_a2r1labs}, this yields a purely deterministic dimension-incremental sparse FFT method. For Algorithm~\ref{algo:sfft_mlfft} and~\ref{algo:sfft_mlfft6}, we still perform probabilistic constructions since we apply Algorithm~\ref{alg:construct_mr1l_I_distinct_primes:mod} and \cite[Algorithm~6]{Kae17} for obtaining reconstructing multiple rank-1 lattices, respectively. We discuss the resulting sample and arithmetic complexities for Algorithm~\ref{algo:sfft_mlfft} and~\ref{algo:sfft_a2r1labs} in Section~\ref{sec:theo:deterministic}.
\end{Remark}

\begin{Remark}\label{rem:function}
In the description of Algorithm~\ref{algo:sfft_mlfft}, \ref{algo:sfft_mlfft6} and~\ref{algo:sfft_a2r1labs}, we mention that the input parameter $p$ is a trigonometric polynomial as a black box.
Alternatively, we may also insert a function $f\in L_1(\T^d)\cap\mathcal{C}(\T^d)$ with infinitely many non-zero Fourier coefficients $\hat{f}_\boldk$.
If these Fourier coefficients decay in a certain way, we may also successfully detect the approximately largest ones by choosing the input parameters appropriately as the numerical results in Section~\ref{sec:numerics:fct} indicate.
\end{Remark}

After presenting the new algorithms, we give its sample and arithmetic complexities.

\subsection{Number of samples and arithmetic complexity}
\label{sec:method:complexity}

We provide a detailed description of the complexities for various steps of Algorithm~\ref{algo:sfft_mlfft}, \ref{algo:sfft_mlfft6} and~\ref{algo:sfft_a2r1labs}. 
We perform the analysis for the case where the search domain~$\Gamma$ is the full grid $\hat{G}_N^d$ of expansion $N\in\N$ and the sparsity parameter $s\gtrsim N$.
For other cases, the resulting sample and arithmetic complexities can be determined analogously. When estimating the complexities, we keep track of the number of detection iterations $r\in\N$ and the theoretical size of this parameter is discussed later in Section~\ref{sec:theoretical}. In practice, as the numerical results in Section~\ref{sec:numerics} suggest, the number of detection iterations~$r$ may be chosen as a constant value in many cases, e.g. $r=1$ for multivariate trigonometric polynomials $p$ in the noiseless case or $\leq 10$ for other considered test settings.

\begin{sloppypar}
In step \ref{enum:reco:incr_general:I2} of Algorithm~\ref{algo:sfft_mlfft} in dimension increment step $t$,
the index sets $I^{(1,\ldots,t-1)}$ and $I^{(t)}$ consist of at most $rs$ and $|\hat{G}_N^1|=2N+1$ many frequencies, respectively.
This yields that the index set $I^{(1,\ldots,t-1)} \times I^{(t)}$ consists of 
$\vert I^{(1,\ldots,t-1)} \times I^{(t)}\vert \leq rs\,(2N+1) \lesssim rsN$
frequency candidates.
The sampling set $\mathcal{X}^{(1,\ldots,t)}$ built in step~\ref{enum:reco:incr_general:I2:X_1_t} of Algorithm~\ref{algo:sfft_mlfft}, where $\mathcal{P}_{(1,\ldots,t)}(\mathcal{X}^{(1,\ldots,t)})$ is a multiple rank-1 lattice returned by applying Algorithm \ref{alg:construct_mr1l_I_distinct_primes:mod} on the index set $J_t:=(I^{(1,\ldots,t-1)} \times I^{(t)})\cap\mathcal{P}_{(1,\ldots,t)}(\Gamma)$, has the size
$$
|\mathcal{X}^{(1,\ldots,t)}| < 2 \tilde{L} \max\{2(rs\,(2N+1)-1), 4\tilde{L} \, \ln\tilde{L} \} \lesssim rsN \log (rsN),
$$
where $\tilde{L}:=\left\lceil 2 \, \ln (rs\,(2N+1)) + 2\,\ln 2\right\rceil$, cf. the proof of Lemma~\ref{lem:reco_mr1l:rep_b} with $c:=2$ and $\gamma:=1/2$.
Moreover, the construction of the sampling set $\mathcal{X}^{(1,\ldots,t)}$ requires
$\mathcal{O}\big(b\,rsN(\log(rsN)+t)\log(rsN)+d\big)$
arithmetic operations.
The inverse multiple rank-1 lattice FFT, cf. Algorithm~\ref{algo:mlfft2:ifft_direct}, in step~\ref{enum:reco:incr_general:I2:p_hat_tilde} requires
$\lesssim\vert\mathcal{X}^{(1,\ldots,t)}\vert\log\vert\mathcal{X}^{(1,\ldots,t)}\vert + \tilde{L} \vert J_t\vert\,(t+\log{|J_t|}) \lesssim rsN \big(\log (rsN) + t\big)\log (rsN)$
arithmetic operations for each detection iteration $i\in\{1,\ldots,r\}$ and each dimension increment step $t\in\{2,\ldots,d\}$.

\end{sloppypar}

In total, this yields
\begin{equation}
\lesssim dr^2sN \log (rsN) + d\,r\,N \lesssim dr^2sN \log (rsN)\label{eq:sample_comp_mr1l}
\end{equation}
many samples and
\begin{equation}
\lesssim d\,r^2 sN  \; (d+\log(rsN)) \, \log (rsN) \lesssim d^2\,r^2 sN \, \log^2 (rsN)\label{eq:arith_comp_mr1l}
\end{equation}
arithmetic operations for Algorithm~\ref{algo:sfft_mlfft} assuming $r \gtrsim b$.
We remark that we do not have any exponential dependence in the dimension $d$, neither for the number of samples nor for the arithmetic complexity,
but only a polynomial dependence of degree $\leq 2$. The sample and arithmetic complexities are distinctly reduced by using Algorithm~\ref{algo:sfft_mlfft} (multiple rank-1 lattices) and compared to Algorithm~\ref{algo:sfft_a2r1labs} or \cite[Algorithm~2]{PoVo14} (single rank-1 lattices), where the latter two algorithms require
$\lesssim dr^3s^2N$ samples and $\lesssim dr^3s^3 + dr^3s^2N\log (rsN)$ arithmetic operations.

In particular, the number of used sampling values as well as the arithmetic complexity is
distinctly improved to an almost linear term with respect to the sparsity parameter $s$
for Algorithm \ref{algo:sfft_mlfft} and multiple rank-1 lattices.

When analyzing the sample and arithmetic complexity of Algorithm~\ref{algo:sfft_mlfft6}, we obtain similar numbers except for additional logarithmic factors.
In total, we require
$$\lesssim dr^2sN \log^2 (rsN) + d\,r\,N \lesssim dr^2sN \log^2 (rsN)$$
many samples with high probability and
$$\lesssim drsN \big(d \log(rsN) + r\log(rsN) + d\big) \log^2 (rsN) \lesssim d^2r^2sN \, \log^3(rsN)$$
arithmetic operations with high probability for Algorithm~\ref{algo:sfft_mlfft6}. We remark that we do not observe these additional logarithmic factors in numerical tests.
Presumably, this is due to the relatively rough estimates of the complexities in~\cite{Kae17}.
Moreover, the numerical results in Section~\ref{sec:numerics} suggest that Algorithm~\ref{algo:sfft_mlfft6} requires even less samples and has smaller runtime than Algorithm~\ref{algo:sfft_mlfft} in practice.

We remark that for obtaining the estimates of the sample and arithmetic complexities in this section, we did not assume that the function under consideration is a multivariate trigonometric polynomial $p$, see also Remark~\ref{rem:function}. In the next section, we will especially address the special case, when the function under consideration is sparse in frequency domain, and we will give possible choices for the input parameters of Algorithm~\ref{algo:sfft_mlfft} and~\ref{algo:sfft_a2r1labs} depending on the targeted success probability.

\FloatBarrier

\section{Theoretical results}\label{sec:theoretical}
In this section, we analyze in detail the sample and arithmetic complexities of Algorithm~\ref{algo:sfft_mlfft} and~\ref{algo:sfft_a2r1labs} applied to multivariate trigonometric polynomials~$p$.
We start with the case of sparse multivariate trigonometric polynomials~$p$ having arbitrary Fourier coefficients~$\hat{p}_\boldk$ in Section~\ref{sec:theoretical:prob}. Afterwards, in Section~\ref{sec:theo:discussion}, we discuss the obtained results.
Later, in Section~\ref{sec:theo:deterministic}, we consider the case where the Fourier coefficients~$\hat{p}_\boldk$ fulfill certain stronger properties, resulting in a modification of Algorithm~\ref{algo:sfft_mlfft} and~\ref{algo:sfft_a2r1labs} with small sample and arithmetic complexities.

\subsection{Multivariate trigonometric polynomial with arbitrary Fourier coefficients}\label{sec:theoretical:prob}

We analyze the computation of (projected) Fourier coefficients in step $t$. We denote
the sampling set $\Lambda=\Lambda(\boldsymbol z_1,M_1,\ldots,\boldsymbol z_{L},M_{L})$ for short and assume that
$\Lambda$ is a reconstructing multiple rank-1 lattice for the frequency index set
$$J_t:=\left(I^{(1, \ldots,t-1)}\times I^{(t)}\right)\cap\mathcal{P}_{(1,\ldots,t)}(\Gamma),$$
$J_t\supset\mathcal{P}_{(1,\ldots,t)}(\operatorname{supp}\hat{p})$, determined by Algorithm \ref{alg:construct_mr1l_I_distinct_primes:mod} in Section \ref{sec:pre:multiple_r1l} or by one of Algorithms 1 to 4 from \cite{Kae17}.
We denote $\boldsymbol{\tilde{x}}:=(x'_{t+1},\ldots,x'_{d})$.
We use Algorithm \ref{algo:mlfft2:ifft_direct} in order to determine the projected Fourier coefficients
$\tilde{\hat{p}}_{(1,\ldots,t),\boldsymbol{k}}^{\Lambda}$, $\boldk\in J_t$.
Due to the construction of the reconstructing multiple rank-1 lattice~$\Lambda$,
for each $\boldk\in J_t$, there exists at least one $\ell\in\{1,\ldots,L\}$ such that
\begin{align}
\boldk\cdot\boldz_\ell\not\equiv\boldk'\cdot\boldz_\ell\imod{M_\ell}
\quad\text{for all $\boldk\neq\boldk'\in J_t$}
\label{eq:r1l_single_k_nonaliasing}
\end{align}
holds and we have
\begin{align*}
\tilde{\hat{p}}_{(1,\ldots,t),\boldsymbol{k}}^{\Lambda(\boldz_\ell,M_\ell)}&:=
\frac{1}{M_{\ell}}\sum_{j=1}^{M_{\ell}-1}
p\left(\frac{j}{M_{\ell}}\boldz_{\ell},\boldsymbol{\tilde{x}}\right)\e^{-2\pi\ii\boldk\cdot\boldz_{\ell}\frac{j}{M_{\ell}}}\\
&=
\sum_{\substack{\boldh\in\operatorname{supp}\hat{p}}}\frac{1}{M_{\ell}}\sum_{j=1}^{M_{\ell}-1}
\hat{p}_{\boldh}\e^{2\pi\ii\boldh\cdot(\boldz_{\ell}\frac{j}{M_{\ell}},\boldsymbol{\tilde{x}})} \, \e^{-2\pi\ii\boldk\cdot\boldz_{\ell}\frac{j}{M_{\ell}}}\\
&=
\sum_{\substack{\boldh\in\operatorname{supp}\hat{p} \\ \left((h_1,\ldots,h_t)-\boldk\right)\cdot \boldz_{\ell}\equiv 0 \imod {M_{\ell}}}}
\hat{p}_{\boldh} \, \e^{2\pi\ii(h_{t+1},\ldots,h_d)\cdot\boldsymbol{\tilde{x}}}\\
&=
\sum_{\substack{\boldh\in\operatorname{supp}\hat{p} \\ (h_1,\ldots,h_t)=\boldk}}
\hat{p}_{(\boldk,h_{t+1},\ldots,h_d)} \, \e^{2\pi\ii(h_{t+1},\ldots,h_d)\cdot\boldsymbol{\tilde{x}}}.
\end{align*}
If there exists more than one $\ell\in\{1,\ldots,L\}$ such that \eqref{eq:r1l_single_k_nonaliasing}
holds, then Algorithm \ref{algo:mlfft2:ifft_direct} computes the average of all these coefficients~$\tilde{\hat{p}}_{(1,\ldots,t),\boldsymbol{k}}^{\Lambda(\boldz_\ell,M_\ell)}$, where
the latter are identical. Consequently, we obtain projected Fourier coefficients
$$
\tilde{\hat{p}}_{(1,\ldots,t),\boldsymbol{k}}^{\Lambda}=
\sum_{\substack{\boldh\in\operatorname{supp}\hat{p} \\ (h_1,\ldots,h_t)=\boldk}}
\hat{p}_{(\boldk,h_{t+1},\ldots,h_d)} \, \e^{2\pi\ii(h_{t+1},\ldots,h_d)\cdot\boldsymbol{\tilde{x}}}.
$$
As in \cite[Section 2.2.2]{PoVo14}, we interpret the projected Fourier coefficients~$\tilde{\hat{p}}_{(1,\ldots,t),\boldsymbol{k}}^{\Lambda}$ as
a trigonometric polynomial $g(\boldsymbol{\tilde{x}}):=\tilde{\hat{p}}_{(1,\ldots,t),\boldsymbol{k}}^{\Lambda}$ with respect to $\boldsymbol{\tilde{x}}\in\T^{d-t}$.
In  Algorithm \ref{algo:sfft_mlfft}, we choose $\boldsymbol{\tilde{x}}\in\T^{d-t}$ uniformly at random.
As an improvement to \cite[Theorem 2.5]{PoVo14}, we estimate the probability 
that a projected Fourier coefficient $\tilde{\hat{p}}_{(1,\ldots,t),\boldsymbol{k}}^{\Lambda}$
is of a certain absolute value, i.e., we estimate the measure of the following set
$$
\left\{\boldsymbol{\tilde{x}}\in\T^{d-t}\colon |g(\boldsymbol{\tilde{x}})|=\left|\sum_{\boldh\in\tilde{I}}\hat{g}_\boldh\,\e^{2\pi\ii\boldh\cdot\boldsymbol{\tilde{x}}}\right| < \delta\right\}.
$$
Similarly, for $\Gamma\supset\operatorname{supp}\hat{p}$ and $t\in\{1,\ldots,d\}$, we have for the projected Fourier coefficient
\begin{align}
 \tilde{\hat{p}}_{t,k_t}
 := &
 \sum\limits_{\ell=0}^{K_t-1} p\left(x_1',\ldots,x_{t-1}',\frac{\ell}{K_t},x_{t+1}',\ldots,x_d'\right) \mathrm{e}^{-2\pi\mathrm{i}\ell k_t/K_t} \nonumber\\
 = &
\sum_{\substack{\boldh\in\operatorname{supp}\hat{p} \\ h_t=k_t}}
\hat{p}_{\boldh} \, \e^{2\pi\ii(h_1,\ldots,h_{t-1},h_{t+1},\ldots,h_d)\cdot\boldsymbol{\tilde{x}}} \quad \text{for } k_t\in\mathcal{P}_t(\Gamma)\label{eq:onedim_aliasing_formula}
\end{align}
in step~\ref{enum:reco:incr_general:I1} and~\ref{enum:reco:incr_general:I2:I_t} of Algorithm~\ref{algo:sfft_mlfft},
and we interpret $\tilde{\hat{p}}_{t,k_t}$ as a trigonometric polynomial $g(\boldsymbol{\tilde{x}}) := \tilde{\hat{p}}_{t,k_t}$
with respect to $\boldsymbol{\tilde{x}}:=(x_1',\ldots,x_{t-1}',x_{t+1}',\ldots,x_d')\in\T^{d-1}$.

\begin{Lemma}\label{lem:prob_g_ineq}
 Let a trigonometric polynomial $g:\T^n\rightarrow\C$, $n\in\N$, $g(\boldx):=\sum_{\boldh\in\tilde{I}} \hat{g}_\boldh \, \mathrm{e}^{2\pi\mathrm{i}\boldh\cdot\boldx}\not\equiv 0$,
 $\tilde{I}\subset\Z^n$, $\vert\tilde{I}\vert < \infty$, be given
 such that the property $\Vert g\vert L^1(\T^n)\Vert > \delta$
for a threshold $\delta>0$  is fulfilled.
 Moreover, let $X_1,\ldots,X_n\in\T$ be independent, identical, uniformly distributed random variables
 and we denote by $\boldsymbol{X}:=(X_1,\ldots,X_n)^\top\in\T^n$ the random vector.
 Then, the probability
 \begin{align*}
  \mathbb{P}(\vert g(\boldsymbol{X})\vert < \delta)
  &\leq
  1-\frac{\|g|L_1(\T^n)\|-\delta}{\|g|L_\infty(\T^n)\|}
  =\frac{\|g|L_\infty(\T^n)\|-\|g|L_1(\T^n)\|+\delta}{\|g|L_\infty(\T^n)\|}
  < 1.
 \intertext{If $\max_{\boldh\in\tilde{I}}\vert\hat{g}_\boldh\vert > \delta$, then} \quad
  \mathbb{P}(\vert g(\boldsymbol{X})\vert < \delta)
  &\leq
  1-\frac{\max_{\boldh\in\tilde{I}}\vert\hat{g}_\boldh\vert-\delta}{\sum_{\boldh\in\tilde{I}}\vert\hat{g}_\boldh\vert}
  =:q
  < 1.
 \end{align*}
 Choosing $r$ random vectors $\boldsymbol{X}_1,\ldots,\boldsymbol{X}_r\in\T^n$ independently, we observe
 \begin{align*}
  \mathbb{P}\left(\bigcap_{i=1}^r\left\{\vert g(\boldsymbol{X}_i)\vert < \delta\right\}\right)\le q^r.
 \end{align*}
\end{Lemma}

\begin{proof}
Let $Y$ be a real valued random variable.
We refer to the lower bound
\begin{align*}
\frac{\mathbb{E}h(Y)-h(\delta)}{\|h(Y)|L_\infty(\R)\|}\le\mathbb{P}(|Y|\ge\delta)
\end{align*}
in \cite[Par. 9.3.A]{Lo77} and apply this for the even and on $[0,\infty)$ nondecreasing function $h(t):=|t|$.
We set $Y=|g(\boldsymbol{X})|$ and we obtain
\begin{align*}
\mathbb{P}(|g(\boldsymbol{X})|\ge\delta)&\ge\frac{\|g|L_1(\T^n)\|-\delta}{\|g|L_\infty(\T^n)\|}
\intertext{and}
\mathbb{P}(|g(\boldsymbol{X})|<\delta)&\le 1-\frac{\|g|L_1(\T^n)\|-\delta}{\|g|L_\infty(\T^n)\|}
 =\frac{\|g|L_\infty(\T^n)\|-\|g|L_1(\T^n)\|+\delta}{\|g|L_\infty(\T^n)\|}.
\end{align*}
 The estimate
 $\vert\hat{g}_\boldh\vert = \vert\int_{\T^n}g(\boldx) \,\mathrm{e}^{-2\pi\mathrm{i}\boldh\cdot\boldx} \mathrm{d}\boldx\vert \leq \int_{\T^n} \vert g(\boldx) \vert \mathrm{d}\boldx = \Vert g\vert L^1(\T^n)\Vert$
 for all $\boldh\in\tilde{I}$ yields
 $\max_{\boldh\in\tilde{I}}\vert\hat{g}_\boldh\vert \leq \Vert g\vert L^1(\T^n)\Vert$.
 Since we have
 $\Vert g\vert L^\infty(\T^n)\Vert = \esssup_{\boldx\in\T^n} \vert g(\boldx)\vert \leq \sum_{\boldh\in\tilde{I}}\vert\hat{g}_\boldh\vert$
 and $\max_{\boldh\in\tilde{I}}\vert\hat{g}_\boldh\vert > \delta$,
 the assertion follows.
\end{proof}

Based on these estimates, we determine bounds on the probability that a projected Fourier coefficient is below an absolute threshold $\delta$.
Since the projected Fourier coefficient may consist only of minimal non-zero Fourier coefficients $\hat{p}_\boldh$, the threshold 
$\delta$ has to be chosen 
$\delta < \min_{\boldh\in\mathrm{supp}\,\hat{p}} \vert \hat{p}_\boldh\vert$ in order to apply the second part of
Lemma~\ref{lem:prob_g_ineq}.

\begin{Corollary}\label{cor:prob_AtC}
 Let a threshold value $\delta>0$, a trigonometric polynomial $p$
 with the property $\min_{\boldh\in\mathrm{supp}\,\hat{p}} \vert \hat{p}_\boldh\vert > \delta$
 and a search space $\Gamma\supset \mathrm{supp}\,\hat{p}$ of finite cardinality be given.
 For fixed $t\in\{1,\ldots,d\}$ and $K_t:=\max(\mathcal{P}_t(\Gamma))-\min(\mathcal{P}_t(\Gamma))+1$, we compute the one-dimensional projected Fourier coefficients for the $t$-th component by
\begin{align*}\tilde{\hat{p}}_{t,k_t}&=\tilde{p}_{t,k_t}(x_1',\ldots,x_{t-1}',x_{t+1}',\ldots,x_d') \\
 &:= \sum\limits_{\ell=0}^{K_t-1} p\left(x_1',\ldots,x_{t-1}',\frac{\ell}{K_t},x_{t+1}',\ldots,x_d'\right) \mathrm{e}^{-2\pi\mathrm{i}\ell k_t/K_t}, \quad k_t\in\mathcal{P}_t(\Gamma),
\end{align*} 
 where the values $x_1',\ldots,x_{t-1}',x_{t+1}',\ldots,x_d'\in\T$ are independently chosen uniformly at random.
 Then, the probability
 \begin{align*}
  \mathbb{P}&(\vert \tilde{\hat{p}}_{t,k_t}\vert < \delta)
\leq
1-
\frac{\|\tilde{\hat{p}}_{t,k_t}|L_1(\T^{d-1})\|-\delta}{\|\tilde{\hat{p}}_{t,k_t}|L_\infty(\T^{d-1})\|} \\
&\le
1-\frac{\max_{\boldh=(h_1,\ldots,h_{t-1},k_t,h_{t+1},\ldots,h_d)\in\operatorname{supp}\hat{p}}|\hat{p}_\boldh|-\delta}{\sum_{\boldh=(h_1,\ldots,h_{t-1},k_t,h_{t+1},\ldots,h_d)\in\operatorname{supp}\hat{p}}|\hat{p}_\boldh|}
=:q_{t,k_t}<1
  \quad
  \text{for } k_t\in\mathcal{P}_t(\mathrm{supp}\,\hat{p})
  \end{align*}
holds
due to Lemma~\ref{lem:prob_g_ineq}.
  Repeating the computation of $\tilde{\hat{p}}_{t,k_t}=\tilde{\hat{p}}_{t,k_t}(\boldsymbol{\tilde{x}})$ for different randomly chosen
  $\boldsymbol{\tilde{x}}=\boldsymbol{\tilde{x}}_1,\ldots,\boldsymbol{\tilde{x}}_r\in\T^{d-1}$, $r\in\N$, we estimate
  $$
  \mathbb{P}\left(\max_{\nu=1,\ldots,r}|\tilde{\hat{p}}_{t,k_t}(\boldsymbol{\tilde{x}}_{\nu})|< \delta\right)
  \le (q_{t,k_t})^r.
  $$
  Applying the union bound yields
  \begin{align}
  \mathbb{P}\left(
  \bigcup_{k_t\in\mathcal{P}_t(\mathrm{supp}\,\hat{p})}
  \left\{\max_{\nu=1,\ldots,r}|\tilde{\hat{p}}_{t,k_t}(\boldsymbol{\tilde{x}}_{\nu})|< \delta\right\}\right)\nonumber
  &\le \sum_{k_t\in\mathcal{P}_t(\mathrm{supp}\,\hat{p})}
  (q_{t,k_t})^r\nonumber \\
  &\le \min\big\{|\operatorname{supp}\hat{p}|,K_t\big\} \left(\max_{k_t\in\mathcal{P}_t(\mathrm{supp}\,\hat{p})}
  q_{t,k_t}\right)^r.\label{eq:cor:prob_AtC}
  \end{align}
\qed
\end{Corollary}

\begin{Corollary}\label{thm:step_t_single_frequency}
Let a threshold value $\delta>0$, a trigonometric polynomial $p$ with the
property $\min_{\boldh\in\operatorname{supp}\hat{p}}|\hat{p}_{\boldh}|>\delta$
and the search space $\Gamma\supset\operatorname{supp}\hat{p}$ of finite cardinality
be given.
For fixed $t\in\{2,\ldots,d-1\}$,
we compute the $t$-dimensional projected Fourier coefficients
$\tilde{\hat{p}}_{(1,\ldots,t),\boldsymbol{k}}^{\Lambda}$,
$\boldk\in\left(I^{(1, \ldots,t-1)}\times I^{(t)}\right)\cap\mathcal{P}_{(1,\ldots,t)}(\Gamma)$, 
in step~\ref{enum:reco:incr_general:I2:p_hat_tilde} of Algorithm~\ref{algo:sfft_mlfft} 
by applying Algorithm \ref{algo:mlfft2:ifft_direct} on $p(\circ,\boldsymbol{\tilde{x}})$, where $\boldsymbol{\tilde{x}}$ is choosen uniformly at random in $\T^{d-t}$.
If $\Lambda(\boldz_1,M_1,\ldots,\boldz_L,M_L)$ is a reconstructing multiple rank-1 lattice
for $\mathcal{P}_{(1,\ldots,t)}(\operatorname{supp}\hat{p})$
determined by Algorithm \ref{alg:construct_mr1l_I_distinct_primes:mod} in Section \ref{sec:pre:multiple_r1l} or by one of \cite[Algorithms 1 to 4]{Kae17}, then
the probability
\begin{align}\nonumber
\mathbb{P}(|\tilde{\hat{p}}_{(1,\ldots,t),\boldsymbol{k}}^{\Lambda}|< \delta)
&\leq \nonumber
1-
\frac{\|\tilde{\hat{p}}_{(1,\ldots,t),\boldsymbol{k}}^{\Lambda}|L_1(\T^{d-t})\|-\delta}{\|\tilde{\hat{p}}_{(1,\ldots,t),\boldsymbol{k}}^{\Lambda}|L_\infty(\T^{d-t})\|}\\
&\le \nonumber
1-\frac{\max_{\boldh=(\boldk,h_{t+1},\ldots,h_d)\in\operatorname{supp}\hat{p}}|\hat{p}_\boldh|-\delta}{\sum_{\boldh=(\boldk,h_{t+1},\ldots,h_d)\in\operatorname{supp}\hat{p}}|\hat{p}_\boldh|}
\\&=:q_{(1,\ldots,t),\boldk}<1 \label{def:q:step_t_single_frequency}
\quad\text{for $\boldk\in \mathcal{P}_{(1,\ldots,t)}(\operatorname{supp}\hat{p})\cap\left(I^{(1, \ldots,t-1)}\times I^{(t)}\right)$}
\end{align}
due to Lemma~\ref{lem:prob_g_ineq}.
Repeating the computation of $\tilde{\hat{p}}_{(1,\ldots,t),\boldsymbol{k}}^{\Lambda}=\tilde{\hat{p}}_{(1,\ldots,t),\boldsymbol{k}}^{\Lambda}(\boldsymbol{\tilde{x}})$ for different randomly chosen
$\boldsymbol{\tilde{x}}=\boldsymbol{\tilde{x}}_1,\ldots,\boldsymbol{\tilde{x}}_r\in\T^{d-t}$, $r\in\N$, we estimate
$$
\mathbb{P}\left(\max_{\nu=1,\ldots,r}|\tilde{\hat{p}}_{(1,\ldots,t),\boldsymbol{k}}^{\Lambda}(\boldsymbol{\tilde{x}}_{\nu})|< \delta\right)
\le (q_{(1,\ldots,t),\boldk})^r.
$$
Applying the union bound yields
\begin{align}
\mathbb{P}&\left(
\bigcup_{\boldk\in\mathcal{P}_{(1,\ldots,t)}(\operatorname{supp}\hat{p})\cap\left(I^{(1, \ldots,t-1)}\times I^{(t)}\right)}
\left\{\max_{\nu=1,\ldots,r}|\tilde{\hat{p}}_{(1,\ldots,t),\boldsymbol{k}}^{\Lambda}(\boldsymbol{\tilde{x}}_{\nu})|< \delta\right\}\right)
\nonumber\\
&\le \sum_{\boldk\in\mathcal{P}_{(1,\ldots,t)}(\operatorname{supp}\hat{p})\cap\left(I^{(1, \ldots,t-1)}\times I^{(t)}\right)}
(q_{(1,\ldots,t),\boldk})^r\nonumber\\
&\le|\operatorname{supp}\hat{p}|\left(\max_{\boldk\in\mathcal{P}_{(1,\ldots,t)}(\operatorname{supp}\hat{p})\cap\left(I^{(1, \ldots,t-1)}\times I^{(t)}\right)}
q_{(1,\ldots,t),\boldk}\right)^r.\label{eq:cor:prob_CtC}
\end{align}
\qed
\end{Corollary}

Analyzing Algorithm \ref{algo:sfft_mlfft} in detail, we
compute the adaptive approximation in $d$ dimension increment steps.
In each of the dimension increment steps, at most three probabilistic sub-steps are performed,
cf.\ step \ref{enum:reco:incr_general:I1}, \ref{enum:reco:incr_general:I2:I_t}, \ref{enum:reco:incr_general:I2:X_1_t}.
If each of these probabilistic sub-steps is successful, we
detect all frequencies from $\operatorname{supp}\hat{p}$ correctly.
We use the union bound to estimate the corresponding probability
\begin{align}\nonumber
\mathbb{P}\left(
\bigcap_{t=1}^d A_t
\cap
\bigcap_{t=2}^d B_t\cap
\bigcap_{t=2}^{d}C_t
\right)
&=1-\mathbb{P}\left(
\bigcup_{t=1}^d A_t^\complement
\cup
\bigcup_{t=2}^d B_t^\complement
\cup
\bigcup_{t=2}^{d}C_t^\complement
\right)\\
&\ge 1- \sum_{t=1}^d \mathbb{P}(A_t^\complement) - \sum_{t=2}^d \mathbb{P}(B_t^\complement) - \sum_{t=2}^{d} \mathbb{P}(C_t^\complement),
\label{eq:success_prob_union_bound}
\end{align}
where the events
\begin{align*}
A_t:=\Big\{&\mathcal{P}_t(\mathrm{supp}\,\hat{p})\subset I^{(t)}\Big\},\\
B_t:=\Big\{&\Lambda(\boldz_1,M_1,\ldots,\boldz_L,M_L) \textnormal{ is a reconstructing multiple rank-1 lattice } \\ &\textnormal{for  }\mathcal{P}_{(1,\ldots,t)}(\operatorname{supp}\hat{p})\Big\},\\
C_t:=\Big\{&\mathcal{P}_{(1,\ldots,t)}(\operatorname{supp}\hat{p})\subset I^{(1, \ldots,t)}\Big\}.
\end{align*}
The probabilities $\mathbb{P}(A_t^\complement)$ and $\mathbb{P}(C_t^\complement)$ of the complements of $A_t$ and $C_t$
were estimated in \eqref{eq:cor:prob_AtC} and \eqref{eq:cor:prob_CtC}, respectively.
An upper bound on $\mathbb{P}(B_t^\complement)$ is obtained via the maximal number~$b$ of multiple rank-1 lattice searches per component.
For the whole algorithm, we allow a failure probability $\varepsilon\in(0,1)$.
We split this up such that each probabilistic sub-step has an equal upper bound on its failure
probability of $\varepsilon/(3d)$.
This strategy allows for estimates of
\begin{itemize}
\item the number of detection iterations~$r$,
\item the maximal number~$b$ of multiple rank-1 lattice searches per component,
\end{itemize}

First, we estimate the required number of detection iterations~$r$.

\begin{Lemma}\label{lem:prob_AtC_CtC}
Let a threshold value $\delta\ge0$, a trigonometric polynomial $p\not\equiv 0$ with the
property $\min_{\boldh\in\operatorname{supp}\hat{p}}|\hat{p}_{\boldh}|\geq 3\delta$
and the search space $\Gamma\supset\operatorname{supp}\hat{p}$ of finite cardinality
be given.
Choosing the number of detection iterations
\begin{align}
r\ge 2|\operatorname{supp}\hat{p}|(\log3+\log d+\log |\operatorname{supp}\hat{p}|-\log\varepsilon)\label{eq:lem:r_estimate}
\end{align}
in Algorithm~\ref{algo:sfft_mlfft}
guarantees that each of the probabilities $\mathbb{P}(A_t^\complement)$ and $\mathbb{P}(C_t^\complement)$
is bounded from above by $\varepsilon/(3d)$.
\end{Lemma}

\begin{proof}
We estimate the probability $\mathbb{P}(C_t^\complement)$ by \eqref{eq:cor:prob_CtC} and we increase $r$
such that
\begin{align*}
\mathbb{P}&\left(
\bigcup_{\boldk\in\mathcal{P}_{(1,\ldots,t)}(\operatorname{supp}\hat{p})\cap\left(I^{(1, \ldots,t-1)}\times I^{(t)}\right)}
\left\{\max_{\nu=1,\ldots,r}|\tilde{\hat{p}}_{(1,\ldots,t),\boldsymbol{k}}^{\Lambda}(\boldsymbol{\tilde{x}}_{\nu})|< \delta\right\}\right)
\nonumber\\
&\le|\operatorname{supp}\hat{p}|\left(\max_{\boldk\in\mathcal{P}_{(1,\ldots,t)}(\operatorname{supp}\hat{p})\cap\left(I^{(1, \ldots,t-1)}\times I^{(t)}\right)}
q_{(1,\ldots,t),\boldk}\right)^r\le \frac{\varepsilon}{3d}.
\end{align*}
is fulfilled. Consequently $r$ has to be bounded from below by
\begin{align}
\frac{\log3+\log d+\log|\operatorname{supp}\hat{p}|-\log\varepsilon}
{-\log\max_{\boldk\in\mathcal{P}_{(1,\ldots,t)}(\operatorname{supp}\hat{p})\cap\left(I^{(1, \ldots,t-1)}\times I^{(t)}\right)}
q_{(1,\ldots,t),\boldk}}.\label{eq:lem:r_CtC_q}
\end{align}
We denote the index set $\tilde{I}:=\mathcal{P}_{(1,\ldots,t)}(\operatorname{supp}\hat{p})\cap\left(I^{(1, \ldots,t-1)}\times I^{(t)}\right)$ and we estimate
\begin{eqnarray*}
&&\left(\max_{\boldk\in\tilde{I}}
q_{(1,\ldots,t),\boldk}\right)^{-1}\\
&\overset{\eqref{def:q:step_t_single_frequency}}{=}&\min_{\boldk\in\tilde{I}}
\frac{\sum_{\boldh=(\boldk,h_{t+1},\ldots,h_d)\in\operatorname{supp}\hat{p}}|\hat{p}_\boldh|}{\sum_{\boldh=(\boldk,h_{t+1},\ldots,h_d)\in\operatorname{supp}\hat{p}}|\hat{p}_\boldh|-\max_{\boldh=(\boldk,h_{t+1},\ldots,h_d)\in\operatorname{supp}\hat{p}}|\hat{p}_\boldh|+\delta}\\
&=&1+\min_{\boldk\in\tilde{I}}
\frac{\max_{\boldh=(\boldk,h_{t+1},\ldots,h_d)\in\operatorname{supp}\hat{p}}|\hat{p}_\boldh|-\delta}{\sum_{\boldh=(\boldk,h_{t+1},\ldots,h_d)\in\operatorname{supp}\hat{p}}|\hat{p}_\boldh|}\\
&\ge& 1+\min_{\boldk\in\tilde{I}}
\frac{2\max_{\boldh=(\boldk,h_{t+1},\ldots,h_d)\in\operatorname{supp}\hat{p}}|\hat{p}_\boldh|}{3\sum_{\boldh=(\boldk,h_{t+1},\ldots,h_d)\in\operatorname{supp}\hat{p}}|\hat{p}_\boldh|}
\ge 1+\frac{2}{3|\operatorname{supp}\hat{p}|}=\frac{3|\operatorname{supp}\hat{p}|+2}{3|\operatorname{supp}\hat{p}|}.
\end{eqnarray*}
Consequently, we obtain
\begin{align*}
&\frac{1}{-\log \left(\max_{\boldk\in\tilde{I}} q_{(1,\ldots,t),\boldk}\right)}
=
\frac{1}{\log \left(\left(\max_{\boldk\in\tilde{I}} q_{(1,\ldots,t),\boldk}\right)^{-1}\right)}
\\
&\le
\frac{1}{\log \left(\frac{3|\operatorname{supp}\hat{p}|+2}{3|\operatorname{supp}\hat{p}|}\right)}
=\frac{1}{\log(3|\operatorname{supp}\hat{p}|+2)-\log(3|\operatorname{supp}\hat{p}|)}
<2|\operatorname{supp}\hat{p}|.
\end{align*}
Choosing $r$ as in \eqref{eq:lem:r_estimate} satisfies the lower bound \eqref{eq:lem:r_CtC_q} and 
$\mathbb{P}(C_t^\complement)\le\varepsilon/(3d)$ is fulfilled.
The bound $\mathbb{P}(A_t^\complement)\le\varepsilon/(3d)$ follows
analogously.
\end{proof}

Next, we deal with the choice of the maximal number~$b$ of multiple rank-1 lattice searches per component.

\begin{Lemma}\label{lem:reco_mr1l:rep_b}
\sloppy
Let a frequency index set $\tilde{I}\subset\Z^d$, $|\tilde{I}|<\infty$, be given.
We apply Algorithm~\ref{alg:construct_mr1l_I_distinct_primes:mod} on $\tilde{I}$ with fixed minimal oversampling factor $c>1$ and fixed failure probability $\gamma\in(0,1)$
in order to try determining a 
reconstructing multiple rank-1 lattice $\Lambda=\Lambda(\boldz_1,M_1,\ldots,\boldz_L,M_L)$.
We repeatedly perform this at most $b\in\N$ many times and stop if $\Lambda$ fulfills properties~\eqref{equ:cond:reco_ml_mlfft2}. 
Then,
the total lattice size $M=|\Lambda(\boldz_1,M_1,\ldots,\boldz_L,M_L)|$ is in
$
\mathcal{O}\big(\max\{|\tilde{I}|,N_{\tilde{I}}\}\,\log|\tilde{I}|\big),
$
and the total number of arithmetic operations is in
$\mathcal{O}\big(b\,\max\{|\tilde{I}|,N_{\tilde{I}}\}(\log|\tilde{I}| + d + \log\log N_{\tilde{I}})\log|\tilde{I}|\big)$,
where $N_{\tilde{I}}:=\max_{j=1,\ldots, d}\{\max_{\zb k\in \tilde{I}}k_j-\min_{\zb l\in \tilde{I}}l_j\}$ is the expansion of $\tilde{I}$.
Moreover,
$$b=\left\lceil
\frac{\log 3+\log d-\log\varepsilon}{|\log\gamma|}
\right\rceil
$$
calls of Algorithm~\ref{alg:construct_mr1l_I_distinct_primes:mod} guarantee
that $\Lambda$ is a reconstructing multiple rank-1 lattice $\Lambda(\boldz_1,M_1,\ldots,\boldz_L,M_L)$
with probability at least $1-\varepsilon/(3d)$, $\varepsilon\in (0,1)$.
\end{Lemma}
\begin{proof}
\sloppy
With probability at least $1-\gamma$, $\Lambda(\boldz_1,M_1,\ldots,\boldz_L,M_L)$ returned by Algorithm~\ref{alg:construct_mr1l_I_distinct_primes:mod} is a reconstructing multiple rank-1 lattice for $\tilde{I}$.
From \eqref{eq:algo_mr1l_construction_sample_compl} and \eqref{eq:algo_mr1l_construction_arithm_compl} for fixed $\gamma$, we immediately obtain the claimed upper bound
for the total lattice size $M$ and the number of arithmetic operations, respectively.
\newline
Repeatedly applying Algorithm~\ref{alg:construct_mr1l_I_distinct_primes:mod} in case of failure up to totally $\tilde{b}$ times results in a success property of at least $1-\gamma^{\tilde{b}}$.
Consequently, ensuring a success property of at least $1-\varepsilon/(3d)$ yields the claimed upper bound $b$ on the number of calls.
\end{proof}

Finally, we combine the obtained results in order to estimate the sample and arithmetic complexity of Algorithm~\ref{algo:sfft_mlfft} for a given success probability $1-\varepsilon$.

\begin{Theorem}\label{thm:r_complexity}
\sloppy
Let a failure probability $\varepsilon\in(0,1)$ be given.
We apply Algorithm \ref{algo:sfft_mlfft} on a
multivariate trigonometric polynomial $p$ using the following parameters.
We choose 
\begin{itemize}
\item the search space $\Gamma$ such that $\operatorname{supp}\hat{p}\subset\Gamma\subset\hat{G}_N^d$,
\item an absolute threshold $0<\delta\le\min_{\boldh\in\operatorname{supp}\hat{p}} |\hat{p}_\boldh|/3$,
\item the sparsity parameter $s\geq |\operatorname{supp}\hat{p}|$,
\item the number of detection iterations $r:=\lceil 2\, |\operatorname{supp}\hat{p}|(\log3+\log d+\log |\operatorname{supp}\hat{p}|-\log\varepsilon) \rceil$,
\item the number of multiple rank-1 lattice searches $b:=\left\lceil
(\log 3+\log d-\log\varepsilon)/\log 2\right\rceil$.
\end{itemize}
Moreover, we assume $|\operatorname{supp}\hat{p}| \gtrsim N$ and $|\operatorname{supp}\hat{p}|\gtrsim d$.
Then, with probability $1-\varepsilon$,
the output of the index set $I$ of Algorithm~\ref{algo:sfft_mlfft} is $\operatorname{supp}\hat{p}$,
the total number of sampling nodes is in $\OO{d\,|\operatorname{supp}\hat{p}|^2N(\log |\operatorname{supp}\hat{p}|)^2|\log\varepsilon|}$,
and the total number of arithmetic operations is in
$\OO{d^2 |\operatorname{supp}\hat{p}|^2N(\log |\operatorname{supp}\hat{p}|)^3|\log\varepsilon|}$.
\newline
When choosing the sparsity parameter \hbox{$s\leq C \, |\operatorname{supp}\hat{p}|$} for a constant \hbox{$C\geq 1$},
we always require
$\lesssim{d\, |\operatorname{supp}\hat{p}|^3 N (\log|\operatorname{supp}\hat{p}|)^3 |\log\varepsilon|\left|\log|\log\varepsilon|\right|}$
many samples and
$\lesssim{d^2 |\operatorname{supp}\hat{p}|^3 N (\log |\operatorname{supp}\hat{p}|)^4 |\log\varepsilon|^2\left|\log|\log\varepsilon|\right|^2}$
arithmetic operations.
\end{Theorem}
\begin{proof}
\sloppy
As discussed above, we use the union bound to estimate
\begin{align*}
\mathbb{P}\left(\operatorname{supp}\hat{p}\subset I\right)
\geq
\mathbb{P}\left(
\bigcap_{t=1}^d A_t
\cap
\bigcap_{t=2}^d B_t\cap
\bigcap_{t=2}^{d}C_t
\right)
&\overset{\eqref{eq:success_prob_union_bound}}{\ge} 1- \sum_{t=1}^d \mathbb{P}(A_t^\complement) - \sum_{t=2}^d \mathbb{P}(B_t^\complement) - \sum_{t=2}^{d} \mathbb{P}(C_t^\complement)
\end{align*}
for the output $I$ of Algorithm~\ref{algo:sfft_mlfft}.
From Lemma~\ref{lem:prob_AtC_CtC}, we obtain
$\mathbb{P}(A_t^\complement)\geq\varepsilon/(3d)$ and $\mathbb{P}(C_t^\complement)\geq\varepsilon/(3d)$.
Moreover, we have $\mathbb{P}(B_t^\complement)\geq\varepsilon/(3d)$ due to Lemma~\ref{lem:reco_mr1l:rep_b}. Altogether, this yields
$\mathbb{P}\left(\operatorname{supp}\hat{p}\subset I\right) \geq 1-\varepsilon$. Next, we discuss the complexities.
\newline
In the general case, the number of frequencies within $I^{(1,\ldots,t)}$ is bounded from above by $r s$, cf. Section~\ref{sec:method:complexity}. We show that a successful application of Algorithm \ref{algo:sfft_mlfft} on the multivariate trigonometric polynomial $p$ does not detect more than $|\operatorname{supp}\hat{p}|$ different frequencies in each dimension increment step.
In particular, we proof by contradiction that $I^{(t)}=\mathcal{P}_{t}(\operatorname{supp}\hat{p})$ and $I^{(1,\ldots,t)}=\mathcal{P}_{(1,\ldots,t)}(\operatorname{supp}\hat{p})$ holds.
\newline
First, we consider the construction of the frequency index sets $I^{(t)}$, $t=1,\ldots,d$, i.e., steps~\ref{enum:reco:incr_general:I1} and \ref{enum:reco:incr_general:I2:I_t}.
We assume that  there exists $k_t\not\in \mathcal{P}_t(\operatorname{supp}\hat{p})$ such that $|\tilde{\hat{p}}_{t,k_t}|\ge \delta>0$,
\begin{eqnarray*}
\tilde{\hat{p}}_{t,k_t}
 & \overset{\eqref{eq:onedim_aliasing_formula}}{:=} &
 \sum\limits_{\ell=0}^{K_t-1} p\left(x_1',\ldots,x_{t-1}',\frac{\ell}{K_t},x_{t+1}',\ldots,x_d'\right) \mathrm{e}^{-2\pi\mathrm{i}\ell k_t/K_t}\\
& = &
\sum_{\substack{\boldh\in\operatorname{supp}\hat{p} \\ h_t\equiv k_t\imod{K_t}}}
 \hat{p}_{\boldh} \, \e^{2\pi\ii(h_1,\ldots,h_{t-1},h_{t+1},\ldots,h_d)\cdot\boldsymbol{\tilde{x}}},
 \end{eqnarray*}
holds, which implies that there exists at least one $h_t\in\Z\setminus\{k_t\}$, such that
$\hat{p}_{\boldh}\neq 0$ and $h_t\equiv k_t\imod{K_t}$.
Since $K_t:=\max(\mathcal{P}_t(\Gamma))-\min(\mathcal{P}_t(\Gamma))+1$
and $k_t\in\mathcal{P}_t(\Gamma)$,
we observe $h_t\not\in\mathcal{P}_{t}(\Gamma)$ and this yields 
$\Gamma\not\supset{\operatorname{supp}\hat{p}}$, which is in contradiction to the choice of $\Gamma$. 
Consequently, the inclusion $I^{(t)}\subset\mathcal{P}_{t}(\operatorname{supp}\hat{p})$ holds and with probability at least $1-\varepsilon/(3d)$, we have the
equality $I^{(t)}=\mathcal{P}_{t}(\operatorname{supp}\hat{p})$.
\newline
Second, in step \ref{enum:reco:incr_general:I2:non-zero} of Algorithm \ref{algo:sfft_mlfft}, the frequency index sets $I^{(1,\ldots,t)}$ are constructed. 
For the index set of frequency candidates $J_t:=(I^{(1,\ldots, t-1)}\times I^{(t)})\cap\mathcal{P}_{(1,\ldots,t)}(\Gamma)$,
we require $J_t \supset \mathcal{P}_{(1,\ldots,t)}(\operatorname{supp}\hat{p})$
and that $\Lambda$ is a reconstructing multiple rank-1 lattice for $J_t$ fulfilling condition~\eqref{equ:cond:reco_ml_mlfft2}.
Then, for each $\boldk\in J_t$ fulfilling $|\tilde{\hat{p}}_{(1,\ldots,t),\boldsymbol{k}}^{\Lambda}| \ge \delta >0$,
we have $\boldk\in\mathcal{P}_{(1,\ldots,t)}(\operatorname{supp}\hat{p})$.
We show this by contradiction and for this,
we assume that there exists a frequency $\boldk\in J_t \setminus \mathcal{P}_{(1,\ldots,t)}(\operatorname{supp}\hat{p})$ with
$|\tilde{\hat{p}}_{(1,\ldots,t),\boldsymbol{k}}^{\Lambda}| \ge \delta >0$. The projected Fourier coefficient $\tilde{\hat{p}}_{(1,\ldots,t),\boldsymbol{k}}^{\Lambda}$ returned by Algorithm \ref{algo:mlfft2:ifft_direct}
is computed as an arithmetic mean of at least one coefficient
\begin{align*}
\tilde{\hat{p}}_{(1,\ldots,t),\boldsymbol{k}}^{\Lambda(\boldz_\ell,M_\ell)}&:=
\frac{1}{M_{\ell}}\sum_{j=1}^{M_{\ell}-1}
p\left(\frac{j}{M_{\ell}}\boldz_{\ell},\boldsymbol{\tilde{x}}\right)\e^{-2\pi\ii\boldk\cdot\boldz_{\ell}\frac{j}{M_{\ell}}}\\
&=
\sum_{\substack{(\boldsymbol{k'},h_{t+1},\ldots,h_d)\in\operatorname{supp}\hat{p} \\ \boldsymbol{k'}\cdot\boldz_\ell\equiv\boldk\cdot\boldz_\ell\imod{M_\ell}}}
\hat{p}_{(\boldsymbol{k'},h_{t+1},\ldots,h_d)} \, \e^{2\pi\ii(h_{t+1},\ldots,h_d)\cdot\boldsymbol{\tilde{x}}}.
\end{align*}
Since $\tilde{\hat{p}}_{(1,\ldots,t),\boldsymbol{k}}^{\Lambda} \neq 0$, there exists at least one $\ell$
such that $\tilde{\hat{p}}_{(1,\ldots,t),\boldsymbol{k}}^{\Lambda(\boldz_\ell,M_\ell)}\neq 0$, and thus, 
there exists at least one $\hat{p}_{(\boldsymbol{k'},h_{t+1},\ldots,h_d)}\neq 0$, where $(\boldsymbol{k'},h_{t+1},\ldots,h_d)\in\operatorname{supp}\hat{p}$ and $\boldsymbol{k'}\cdot\boldz_\ell\equiv\boldk\cdot\boldz_\ell\imod{M_\ell}$.
On the one hand, due to the assumption $\boldk\not\in \mathcal{P}_{(1,\ldots,t)}(\operatorname{supp}\hat{p})$, we have
$\boldsymbol{k'}\neq\boldk$. On the other hand,
 $\Lambda$ is a reconstructing multiple rank-1 lattice constructed by Algorithm \ref{alg:construct_mr1l_I_distinct_primes:mod}
 fulfilling the properties \eqref{equ:cond:reco_ml_mlfft2}. Consequently, we obtain
 $\boldsymbol{k'}\not\in J_t$, which is in contradiction to the requirement
$J_t \supset \mathcal{P}_{(1,\ldots,t)}(\operatorname{supp}\hat{p})$.
\newline
If all events $A_1$, and $A_t$, $B_t$, $C_t$ for $t\in\{2,\ldots,d\}$ occur, the equalities $I^{(1)}=\mathcal{P}_{1}(\operatorname{supp}\hat{p})$, $I^{(t)}=\mathcal{P}_{t}(\operatorname{supp}\hat{p})$,
and $I^{(1,\ldots,t)}=\mathcal{P}_{(1,\ldots,t)}(\operatorname{supp}\hat{p})$ hold.
In particular, this means that the cardinalities are bounded by $|I^{(1,\ldots,t)}|\le|\operatorname{supp}\hat{p}|$.
This leads to a distinct reduction in the upper bounds of the complexities in \eqref{eq:sample_comp_mr1l} and \eqref{eq:arith_comp_mr1l}.
In more detail, the factors $r\,s$ can be replaced by $|\operatorname{supp}\hat{p}|$, which leads to
$\OO{dr|\operatorname{supp}\hat{p}|N \log (|\operatorname{supp}\hat{p}|\,N)}$ samples and
$\OO{d^2\,r |\operatorname{supp}\hat{p}| N \, \log^2 (|\operatorname{supp}\hat{p}|\,N)}$ arithmetic operations.
\newline
The latter discussed events occur with probability $1-\varepsilon$ as mentioned above.
Taking the assumptions $|\operatorname{supp}\hat{p}| \gtrsim N$ and $|\operatorname{supp}\hat{p}|\gtrsim d$
as well as the choice of the parameters~$r$ and~$b$ into account, we obtain
with probability at least $1-\varepsilon$
that the number of samples is in
$\OO{d\,|\operatorname{supp}\hat{p}|^2N \log^2 |\operatorname{supp}\hat{p}|\,|\log\varepsilon|}$ and
that the number of arithmetic operations is in
$\OO{d^2\,|\operatorname{supp}\hat{p}|^2 N \, \log^3 |\operatorname{supp}\hat{p}|\,|\log\varepsilon|}$.
\newline
The worst case number of samples and arithmetic complexities follow from the parameter choices and the discussions in Section~\ref{sec:method:complexity}.
\end{proof}

Similarly, we obtain estimates for the sample and arithmetic complexity of Algorithm~\ref{algo:sfft_a2r1labs}.

\begin{Theorem}\label{cor:r_complexity:single_r1l}
\sloppy
Let a failure probability $\varepsilon\in(0,1)$ be given.
We apply Algorithm \ref{algo:sfft_a2r1labs} on a
trigonometric polynomial $p$ using the following parameters.
We choose 
\begin{itemize}
\item the search space $\Gamma$ such that $\operatorname{supp}\hat{p}\subset\Gamma\subset\hat{G}_N^d$,
\item an absolute threshold $0<\delta\le\min_{\boldh\in\operatorname{supp}\hat{p}} |\hat{p}_\boldh|/3$,
\item the sparsity parameter $s\geq |\operatorname{supp}\hat{p}|$,
\item the number of detection iterations $r:=\lceil 2|\operatorname{supp}\hat{p}|(\log2+\log d+\log |\operatorname{supp}\hat{p}|-\log\varepsilon) \rceil$.
\end{itemize}
Moreover, we assume $|\operatorname{supp}\hat{p}| \gtrsim N$ and $|\operatorname{supp}\hat{p}|\gtrsim d$.
Then, with probability $1-\varepsilon$,
the output~$I$ of Algorithm~\ref{algo:sfft_mlfft} is $\operatorname{supp}\hat{p}$,
the total number of sampling nodes is in $\OO{d \,|\operatorname{supp}\hat{p}|^3N \, (\log|\operatorname{supp}\hat{p}|) |\log\varepsilon|}$,
and the total number of arithmetic operations is in
$\OO{d\,|\operatorname{supp}\hat{p}|^3 N \,(\log|\operatorname{supp}\hat{p}|)^2 |\log\varepsilon|}$.
\newline
When choosing the sparsity parameter \hbox{$s\leq C \, |\operatorname{supp}\hat{p}|$} for a constant \hbox{$C\geq 1$},
we always require
$\lesssim d\, |\operatorname{supp}\hat{p}|^5 N (\log|\operatorname{supp}\hat{p}|)^3 |\log\varepsilon|^3$
many samples and
$\lesssim{d\, |\operatorname{supp}\hat{p}|^6\,(\log |\operatorname{supp}\hat{p}|)^3\,|\log\varepsilon|^3 + d\, |\operatorname{supp}\hat{p}|^5 \, N (\log |\operatorname{supp}\hat{p}|)^4 |\log\varepsilon|^3\left|\log|\log\varepsilon|\right|}$
arithmetic operations.
\end{Theorem}
\begin{proof}
We proceed similarly as in the proof of Theorem~\ref{thm:r_complexity}. The main difference is that we use reconstructing single rank-1 lattices, which are constructed by a deterministic algorithm.
Consequently, we have 
\begin{align*}
\mathbb{P}\left(\operatorname{supp}\hat{p}\subset I\right)
\geq
\mathbb{P}\left(
\bigcap_{t=1}^d A_t
\cap
\bigcap_{t=2}^{d}C_t
\right)
&\overset{\eqref{eq:success_prob_union_bound}}{\ge} 1- \sum_{t=1}^d \mathbb{P}(A_t^\complement) - \sum_{t=2}^{d} \mathbb{P}(C_t^\complement)
\end{align*}
for the output $I$ of Algorithm~\ref{algo:sfft_mlfft}. For each probability $\mathbb{P}(A_t^\complement)$ and $\mathbb{P}(C_t^\complement)$, we ensure an lower bound of $\varepsilon/(2d)$.
Consequently, we have a slightly smaller choice for the number of detection iterations $r$ in Lemma~\ref{lem:prob_AtC_CtC}.
\newline
For the number of samples, the dominating terms are the sizes $M^{(t)}$ of the single rank-1 lattices in each dimension increment step $t\in\{2,\ldots,d\}$.
If all of the previous steps succeeded, we have $M^{(t)}\lesssim |\operatorname{supp}\hat{p}|^2 N$ by \cite[Corollary~2.3]{PoVo14} as a consequence of \cite{Kae2013}.
In total, this means
$$\lesssim d\,r\,|\operatorname{supp}\hat{p}|^2 N \lesssim d\,|\operatorname{supp}\hat{p}|^3N(\log |\operatorname{supp}\hat{p}|) \,|\log\varepsilon|$$
many samples with probability at least $1-\varepsilon$. \newline
Building the reconstructing single rank-1 lattices requires $\OO{d\,|\operatorname{supp}\hat{p}|^3}$ arithmetic operations in total and computing the rank-1 lattice FFTs
$\OO{d\,r\,|\operatorname{supp}\hat{p}|^2 N \log |\operatorname{supp}\hat{p}|}$ many in total with probability at least $1-\varepsilon$, which yields the claimed arithmetic complexity $\OO{d\,|\operatorname{supp}\hat{p}|^3 N \,(\log|\operatorname{supp}\hat{p}|)^2 |\log\varepsilon|}$.
\newline
Since Algorithm \ref{algo:sfft_a2r1labs} requires
$\lesssim dr^3s^2N$ samples and $\lesssim dr^3s^3 + dr^3s^2N\log (rsN)$ arithmetic operations, cf. \cite[Section 2.2.3]{PoVo14},
the worst case complexities follow due to the parameter choices of $r$ and $s$.
\end{proof}

\subsection{Discussion}\label{sec:theo:discussion}
In Section~\ref{sec:theoretical:prob}, we obtained detailed theoretical estimates for the sample and arithmetic complexities of Algorithm~\ref{algo:sfft_mlfft} and~\ref{algo:sfft_a2r1labs} when applied on multivariate trigonometric polynomials~$p$.
Starting point of the theoretical results was Lemma~\ref{lem:prob_g_ineq}, which estimates the probability that the absolute value of a multivariate trigonometric polynomial is below or equal to a certain threshold~$\delta$ for a uniformly random sample. Afterwards, this result was used in Corollary~\ref{cor:prob_AtC} and~\ref{thm:step_t_single_frequency} to bound the probability that the absolute value of a single aliased Fourier coefficient in Algorithm~\ref{algo:sfft_mlfft} or~\ref{algo:sfft_a2r1labs} is below or equal to a given absolute threshold~$\delta$, which would mean that the corresponding frequency component is not detected. Using a certain number of detection iterations $r$, repeated sampling can be used to decrease the overall failure probability and increase the overall probability of successfully detecting all non-zero Fourier coefficients $\hat{p}_\boldk\neq 0$ of the multivariate trigonometric polynomial~$p$ under consideration. Based on the previously mentioned results, the required number of detection iterations~$r$ was estimated by
$$
r\ge 2|\operatorname{supp}\hat{p}|(\log3+\log d+\log |\operatorname{supp}\hat{p}|-\log\varepsilon)
$$
in Lemma~\ref{lem:prob_AtC_CtC}
for a given failure probability $\leq\varepsilon\in (0,1)$.
This lower bound for $r$ contains the number $|\operatorname{supp}\hat{p}|$ of non-zero Fourier coefficients as a linear factor which is probably overestimated.
Moreover, even if we would decrease the absolute threshold $\delta$ to almost zero, the currently used proof technique would still result in a lower bound which contains $|\operatorname{supp}\hat{p}|$ as a linear factor.
However, this does not match our expectations and the observed numerical results in Section~\ref{sec:numerics}.

We remark that Lemma \ref{lem:prob_g_ineq} still leads to a distinct improvement
compared to the previous results based on \cite[Lemma~2.4]{PoVo14}, which would yield an even larger bound on $r$ containing the factor $|\operatorname{supp}\hat{p}|^2$ in Section~\ref{sec:theoretical:prob}.

In order to obtain a smaller bound on $r$ and, consequently, lower sample and arithmetic complexities, an improved version of Lemma~\ref{lem:prob_g_ineq} would be required. For instance, an answer to the following more general question for trigonometric polynomials would immediately give such an improvement. Assume that a trigonometric polynomial $g\colon\T^n\rightarrow\C$, $n\in\N$, is given
and let
$$
q:=\int_{\{\boldx\in\T^n\colon |g(\boldx)|\leq \delta\}} \mathrm{d}\,\boldx
= \left|\left\{\boldsymbol{\tilde{x}}\in\T^n\colon |g(\boldsymbol{\tilde{x}})|=\left|\sum_{\boldk\in\operatorname{supp}\hat{g}}\hat{g}_\boldk\,\e^{2\pi\ii\boldh\cdot\boldsymbol{\tilde{x}}}\right|\le \delta\right\}\right|.
$$
How is it possible to show
$$
\frac{1}{\log{\frac{1}{q}}} \underset{\text{\large $\sim$}}{\ll} |\operatorname{supp}\hat{g}|
$$
and which conditions on the Fourier coefficients $\hat{g}_\boldk$ are sufficient.
Of course, these conditions should be as weak as possible.

As already mentioned in \cite{PoVo14}, assuming a special property of the Fourier coefficients $\hat{p}_\boldk$ of a multivariate trigonometric polynomial $p$ allows for choosing the number of detection iterations $r:=1$ as we discuss in the next section. However, for a more general case, we currently do not know how to obtain improved bounds on $r$.

\subsection[Multivariate trigonometric polynomial with special Fourier coefficients allowing for r=1 detection iteration]{Multivariate trigonometric polynomial with special Fourier coefficients allowing for $r=1$ detection iteration}\label{sec:theo:deterministic}
If we assume that the function under consideration is a multivariate trigonometric polynomial~$p$ with the additional property that the signs of the real part $\mathrm{Re}(\hat{p}_\boldk)$ of all non-zero Fourier coefficients $\hat{p}_\boldk$ are identical $\neq 0$ or that the signs of the imaginary part $\mathrm{Im}(\hat{p}_\boldk)$ are identical $\neq 0$, then
we obtain a deterministic version of Algorithm~\ref{algo:sfft_a2r1labs} as explained in Remark~\ref{rem:deterministic}. This strategy requires $\OO{d\,|\operatorname{supp}\hat{p}|^2 N}$ samples and
$\OO{d\,|\operatorname{supp}\hat{p}|^3 + d\,|\operatorname{supp}\hat{p}|^2 N}$ arithmetic operations.

\begin{sloppypar}
In this setting, Algorithm~\ref{algo:sfft_mlfft} is not purely deterministic due to the probabilistic search for the reconstructing multiple rank-1 lattices.
From the proof of Theorem~\ref{thm:r_complexity}, we obtain that Algorithm~\ref{algo:sfft_mlfft} detects all non-zero Fourier coefficients with probability
at least $1-\varepsilon$ when setting the number of multiple rank-1 lattice searches parameter $b:=\left\lceil (\log d-\log\varepsilon)/\log 2\right\rceil$
and the sparsity $s\geq|\operatorname{supp}\hat{p}|$.
This requires $\mathcal{O}\big(d\,|\operatorname{supp}\hat{p}|N \log (|\operatorname{supp}\hat{p}|\,N)\big)$ samples and $\mathcal{O}\big(d^2\, |\operatorname{supp}\hat{p}|N \log^2(|\operatorname{supp}\hat{p}|\,N) \,(\log (d) - \log\varepsilon)\big)$ arithmetic operations with probability at least $1-\varepsilon$. Choosing the sparsity parameter $s\leq C\, |\operatorname{supp}\hat{p}|$ for a constant $C\geq 1$ guarantees
these complexities.
\end{sloppypar}

\section{Numerical results}
\label{sec:numerics}

In this section, we present several numerical results, which empirically confirm the high effectiveness and robustness of the proposed method.
In Section~\ref{sec:numerics:rand_poly}, we start with the reconstruction of random multivariate trigonometric polynomials $p$ in up to $d=30$ spatial dimensions, where the sampling values are not perturbed by noise.
Additionally, we consider the case, where the sampling values are perturbed by additive Gaussian noise in Section~\ref{sec:numerics:noise}.
Moreover, in Section~\ref{sec:numerics:fct}, we successfully determine the approximately largest Fourier coefficients and the corresponding frequencies of a 10-dimensional test function.

In all numerical tests, we compare the new algorithms based on multiple rank-1 lattice sampling with
Algorithm~\ref{algo:sfft_a2r1labs}, which uses single reconstructing rank-1 lattices.
On the one hand, we consider Algorithm~\ref{algo:sfft_mlfft} to confirm theoretical results from Section~\ref{sec:theoretical}.
On the other hand, all numerical test settings will be treated by
Algorithm~\ref{algo:sfft_mlfft6}, for which we currently do not have an extensive theoretical framework.
The main difference in both algorithms is the construction of 
the multiple rank-1 lattices.
The construction method used in Algorithm~\ref{algo:sfft_mlfft6} does not allow for the application of
the proof techniques we developed in Section~\ref{sec:theoretical} and may result in
higher theoretical upper bounds on the sampling and arithmetic complexities, see also the discussion in \cite[Section 4.1]{Kae17}.
However, for the numerical test settings in \cite{Kae17},
the number of sampling nodes was distinctly lower, when using this construction method,
compared to the construction method used in Algorithm~\ref{algo:sfft_mlfft}.
Since we observe a similar preferable behavior when using Algorithm~\ref{algo:sfft_mlfft6},
we also present the corresponding results in detail.

We implemented the methods and numerical tests in MATLAB. All tests were run using MATLAB R2015b and IEEE~754 double precision arithmetic.
The time measurements were performed on a computer with Intel Xeon 5160 CPU (3~GHz) and 64~GB RAM using 1 thread.

In general, we observe that the total number of samples is distinctly lower when using multiple rank-1 lattices instead of single rank-1 lattices.

We remark that when performing the numerical tests, we implicitly assumed that we have at least some basic knowledge about the function under consideration, like e.g. a reasonable search domain~$\Gamma$ or estimates of the minimal absolute value of the Fourier coefficients $\hat{p}_\boldk$, the sparsity $\vert\mathrm{supp}\,\hat{p}\vert$, the noise level, etc.

\subsection{Random sparse trigonometric polynomial}
\label{sec:numerics:rand_poly}

As in \cite[Section 3.1]{PoVo14}, we construct random multivariate trigonometric polynomials $p$
with frequencies supported within the cube $\hat{G}_{N}^d=[-N,N]^d\cap\Z^d$.
For this, we choose $\vert\mathrm{supp}\,\hat{p}\vert$ frequencies $\boldk\in\hat{G}_{N}^d$ uniformly at random
and corresponding Fourier coefficients $\hat{p}_\boldk\in[-1,1)+[-1,1)\mathrm{i}$, $\vert\hat{p}_\boldk\vert\geq 10^{-6}$, $\boldk\in I=\mathrm{supp}\,\hat{p}$.
For the reconstruction of the trigonometric polynomials $p$, we only assume that the search domain $\Gamma:=\hat{G}_{N}^d\supset\mathrm{supp}\,\hat{p}$.
Additionally, we do not truncate the frequency index sets of detected frequencies $I^{(1,\ldots,t)}$, $t\in\{2,\ldots,d\}$,
i.e., we set the sparsity parameter $s:=\vert\Gamma\vert$. Alternatively, one may choose $s:=|\mathrm{supp}\,\hat{p}|$ and this yields identical results for the considered random sparse trigonometric polynomials.

\begin{Example}\label{numerics:rand_poly:samples_errors}
(Number of samples and errors).
We set the expansion $N:=32$.
Now, we compare the results of Algorithm~\ref{algo:sfft_a2r1labs}, which is a modification of \cite[Algorithm~2]{PoVo14} with absolute thresholding and which uses reconstructing single rank-1 lattices, with Algorithm~\ref{algo:sfft_mlfft} and~\ref{algo:sfft_mlfft6}, which use multiple rank-1 lattices as sampling sets.
We set the absolute threshold parameter $\delta:=10^{-12}$ and the number of detection iterations $r:=1$.
For sparsity $\vert\mathrm{supp}\,\hat{p}\vert:=$ 1\,000 and 10\,000, we run tests for dimension $d\in\{5,10,15,20,25,30\}$.
All tests are repeated $10$ times with newly chosen frequencies $\boldk\in \Gamma$ and Fourier coefficients $\hat{p}_\boldk\in \mathbb C$.
Then, for the 10 repetitions, we determine the maximum of the total number of samples and the maximum of the relative $\ell_2$-errors of the Fourier coefficients.
For Algorithm~\ref{algo:sfft_a2r1labs} using reconstructing single rank-1 lattices, the results are shown in columns 3--4 of Table~\ref{table:sfft:dimincr:per:numerics:rand_poly:alg2_mlfft4}.
Moreover, the results for Algorithm~\ref{algo:sfft_mlfft} and~\ref{algo:sfft_mlfft6} using multiple rank-1 lattices are presented in columns 5--6 and 7--8, respectively.
For multiple rank-1 lattices, we also run the tests for sparsity $\vert\mathrm{supp}\,\hat{p}\vert:=100\,000$.
We observe that in the considered cases, the obtained relative $\ell_2$-errors are near machine precision. Moreover, we require distinctly less samples than if we had used a $d$-dimensional FFT on a full grid, which would require $|\hat{\Gamma}_{32}^5| = 1\,160\,290\,625$ already in the 5-dimensional case.
Additionally, for sparsity $\vert\mathrm{supp}\,\hat{p}\vert:=10\,000$, Algorithm~\ref{algo:sfft_mlfft} using multiple rank-1 lattices required only about 1/9 of the samples compared to Algorithm~\ref{algo:sfft_a2r1labs} using single rank-1 lattices.
We observe that the number of samples may be even further reduced by using Algorithm~\ref{algo:sfft_mlfft6}, which required only 1/84 to 1/41 of the samples compared to single rank-1 lattices for sparsity $\vert\mathrm{supp}\,\hat{p}\vert:=10\,000$.
For sparsity $\vert\mathrm{supp}\,\hat{p}\vert:=100\,000$, we were not able to apply Algorithm~\ref{algo:sfft_a2r1labs}, which uses single rank-1 lattices, since this would have required too many samples and extremely long runtimes.
\newline
Instead of choosing the absolute threshold parameter $\delta:=10^{-12}$, we also repeated the tests for a larger absolute threshold of 
$\delta:=10^{-7}\leq\min_{\boldk\in\mathrm{supp}\,\hat{p}}\vert\hat{p}_\boldk|/10$ and we obtained almost identical results as in Table~\ref{table:sfft:dimincr:per:numerics:rand_poly:alg2_mlfft4}.
\newline
We remark that we also ran the numerical tests applying a modified version of Algorithm~\ref{algo:sfft_mlfft}, which uses \cite[Algorithm~4]{Kae17} instead of Algorithm~\ref{alg:construct_mr1l_I_distinct_primes:mod} for determining reconstructing multiple rank-1 lattices. Then, we obtained results comparable to those of Algorithm~\ref{algo:sfft_mlfft} in Table~\ref{table:sfft:dimincr:per:numerics:rand_poly:alg2_mlfft4} requiring slightly less samples, i.e., between 2 and 10 percent less in most cases.
\end{Example}

\begin{table}[htb]
\centering
\begin{small}
\begin{tabular}{|r|r||S[table-format=10.0]|r||S[table-format=10.0]|r||S[table-format=9.0]|r|}
\hline
& & \multicolumn{2}{c||}{Algorithm~\ref{algo:sfft_a2r1labs} using}& \multicolumn{2}{c||}{Algorithm~\ref{algo:sfft_mlfft} using} & \multicolumn{2}{c|}{Algorithm~\ref{algo:sfft_mlfft6} using} \\
& & \multicolumn{2}{c||}{single rank-1 lattices}& \multicolumn{2}{c||}{multiple rank-1 lattices} & \multicolumn{2}{c|}{multiple rank-1 lattices} \\
& & \multicolumn{1}{r|}{max.\ total} & max.\ rel. & \multicolumn{1}{r|}{max.\ total} & max.\ rel. & \multicolumn{1}{r|}{max.\ total} & max.\ rel. \\
$d$\rule[0.1em]{0em}{1em} & $\vert\mathrm{supp}\,\hat{p}\vert$
& \multicolumn{1}{r|}{\#samples} & $\ell_2$-error & \multicolumn{1}{r|}{\#samples} & $\ell_2$-error & \multicolumn{1}{r|}{\#samples} & $\ell_2$-error \\
\hline
\hline
5\rule[0.55em]{0em}{0.55em} & 1\,000 & 6260605 & 8.0e-16 &    4525799 & 5.3e-16 &     581881 & 5.8e-16 \\
10 & 1\,000 & 20848685 & 5.4e-16 &   12115199 & 5.3e-16 &    1589349 & 5.1e-16 \\
15 & 1\,000 & 35937525 & 7.2e-16 &   18827483 & 5.3e-16 &    2599029 & 5.2e-16 \\
20 & 1\,000 & 52361205 & 6.9e-16 &   27030181 & 5.3e-16 &    3609753 & 5.2e-16 \\
25 & 1\,000 & 67164695 & 5.0e-16 &   33666561 & 5.3e-16 &    4621205 & 5.1e-16 \\
30 & 1\,000 & 80660385 & 5.7e-16 &   40085363 & 5.3e-16 &    5644059 & 5.1e-16 \\
\hline
5\rule[0.55em]{0em}{0.55em} & 10\,000 & 190618285 & 8.6e-16 &   40854167 & 3.5e-16 &    4648335 & 4.0e-16 \\
10 & 10\,000 & 1081274675 & 6.3e-16 &   129929699 & 3.4e-16 &   15186447 & 3.8e-16 \\
15 & 10\,000 & 1969412575 & 1.3e-15 &  225282365 & 3.5e-16 &   25662189 & 3.8e-16 \\
20 & 10\,000 & 2935663575 & 5.3e-16 & 316957693 & 3.4e-16 &   36161887 & 3.8e-16 \\
25 & 10\,000 & 3837073825 & 7.1e-16 & 409305929 & 3.4e-16 &   46681103 & 3.8e-16 \\
30 & 10\,000 & 4771398905 & 1.3e-15 & 492181865 & 3.5e-16 &   57203659 & 3.8e-16 \\
\hline
   5\rule[0.55em]{0em}{0.55em} & 100\,000 & \multicolumn{1}{c|}{--} & \multicolumn{1}{c||}{--} & 419479177 & 2.2e-16 &   33428113 & 3.0e-16 \\ 
  10 & 100\,000 & \multicolumn{1}{c|}{--} & \multicolumn{1}{c||}{--} & 1445678877 & 2.1e-16 &  143681689 & 2.1e-16 \\ 
  15 & 100\,000 & \multicolumn{1}{c|}{--} & \multicolumn{1}{c||}{--} & 2572013451 & 2.1e-16 &  250232085 & 2.1e-16 \\ 
  20 & 100\,000 & \multicolumn{1}{c|}{--} & \multicolumn{1}{c||}{--} & 3623019125 & 2.0e-16 &  356857499 & 2.1e-16 \\ 
  25 & 100\,000 & \multicolumn{1}{c|}{--} & \multicolumn{1}{c||}{--} & 4678498335 & 2.1e-16 &  463174925 & 2.9e-16 \\ 
  30 & 100\,000 & \multicolumn{1}{c|}{--} & \multicolumn{1}{c||}{--} & 5680886471 & 2.0e-16 &  569711277 & 2.1e-16 \\ 
\hline
\end{tabular}
\end{small}
\caption{Results for random sparse trigonometric polynomials using Algorithm~\ref{algo:sfft_a2r1labs}, \ref{algo:sfft_mlfft}, \ref{algo:sfft_mlfft6},
when considering frequencies within the search domain $\Gamma=\hat{G}_{32}^d$.}
\label{table:sfft:dimincr:per:numerics:rand_poly:alg2_mlfft4}
\end{table}

From the previous example, we observe that for fixed dimension $d$ and increasing sparsity $\vert\mathrm{supp}\,\hat{p}\vert$,
the number of samples increases only moderately for Algorithm~\ref{algo:sfft_mlfft} and~\ref{algo:sfft_mlfft6}.
In the next example, we investigate the behavior in more detail and additionally measure the runtimes using only one CPU thread.

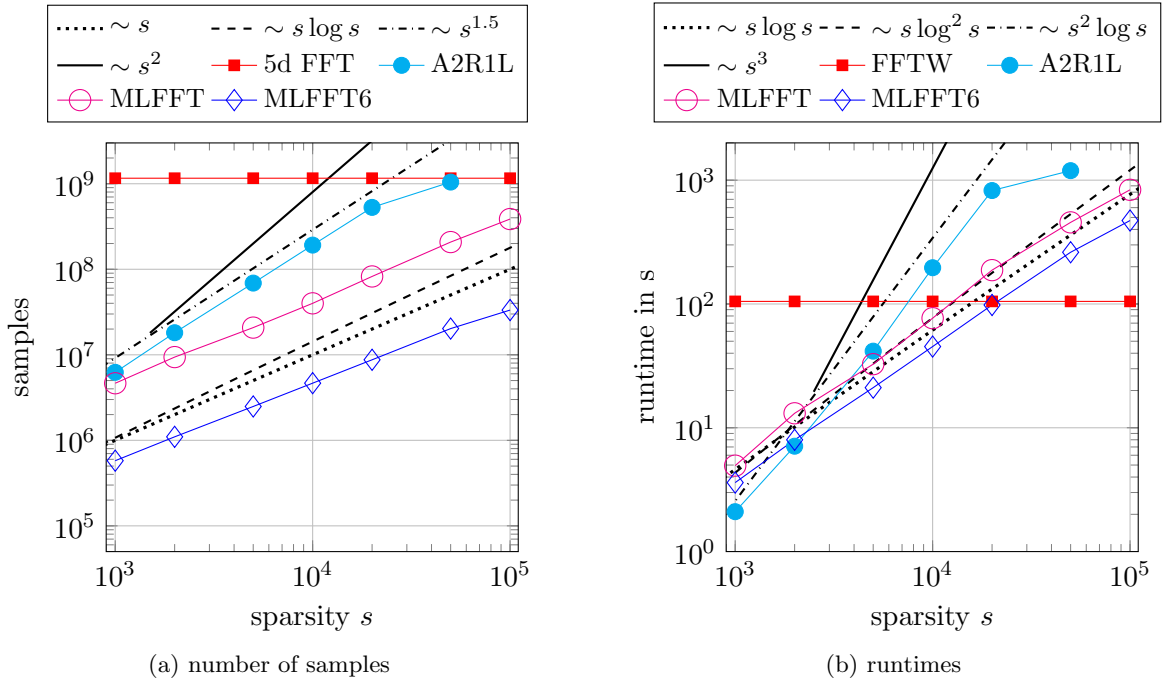
\begin{figure}[htb]
\subfloat[number of samples]{\label{figure:numerics:rand_poly:samples_runtimes:r1:samples}
\begin{tikzpicture}
  \begin{loglogaxis}[enlargelimits=false,xmin=900,xmax=110000,ymin=50000,ymax=3e9,ytick={100,1000,1e4,1e5,1e6,1e7,1e8,1e9},height=0.45\textwidth, width=0.45\textwidth, grid=major, xlabel={sparsity $s$}, ylabel={samples},
   legend style={at={(0.44,1.05)}, anchor=south,legend columns=3,legend cell align=left, font=\small },
  ]
\addplot[black,domain=1e1:1e6,samples=100,dotted,very thick] {x*1000};
\addlegendentry{$\sim s$}
\addplot[black,domain=1e1:1e6,samples=100,dashed,thick] {x*ln(x)*155};
\addlegendentry{$\sim s\log s$}
\addplot[black,domain=1e1:1e6,samples=100,dashdotted,thick] {290*x^1.5};
\addlegendentry{$\sim s^{1.5}$}
\addplot[black,domain=1.5e3:1e5,samples=100,solid,thick] {8*x^2};
\addlegendentry{$\sim s^2$}
\addplot[red,mark=square*,mark size=2] coordinates {
 (1000,1160290625) (2000,1160290625) (5000,1160290625) (10000,1160290625) (20000,1160290625) (50000,1160290625) (100000,1160290625)
};
\addlegendentry{5d FFT}
\addplot[cyan,mark=*,mark size=3] coordinates {
 (1000,6237595) (2000,18121415) (5000,68928405) (10000,191692475) (20000,529037925) (50000,1045516095)
};
\addlegendentry{A2R1L}
\addplot[magenta,mark=o,mark size=4] coordinates {
 (1000,4647753) (2000,9401115) (5000,20796159) (10000,39971549) (20000,82431095) (50000,207916371) (100000,386639121)
};
\addlegendentry{MLFFT}
\addplot[blue,mark=diamond,mark size=4] coordinates {
 (1000,577983) (2000,1098507) (5000,2488017) (10000,4660439) (20000,8740939) (50000,20247793) (100000,33332473)
};
\addlegendentry{MLFFT6}
  \end{loglogaxis}
\end{tikzpicture}
}
\hfill
\subfloat[runtimes]{\label{figure:numerics:rand_poly:samples_runtimes:r1:runtimes}
\begin{tikzpicture}
  \begin{loglogaxis}[enlargelimits=false,xmin=900,xmax=110000,ymin=1,ymax=2000,ytick={1,10,100,1000},height=0.45\textwidth, width=0.45\textwidth, grid=major, xlabel={sparsity $s$}, ylabel={runtime in s},
   legend style={at={(0.44,1.05)}, anchor=south,legend columns=3,legend cell align=left, font=\small },
  ]
\addplot[black,domain=1e1:1e6,samples=100,dotted,very thick] {x*ln(x)/1500};
\addlegendentry{$\sim s\log s$}
\addplot[black,domain=1e3:1e6,samples=100,dashed,thick] {x*(ln(x))^2/11000};
\addlegendentry{$\sim s\log^2 s$}
\addplot[black,domain=1e1:1e5,samples=100,dashdotted,thick] {x^2*(ln(x))/2700000};
\addlegendentry{$\sim s^2\log s$}
\addplot[black,domain=2.5e3:1e5,samples=100,solid,thick] {x^3/800000000};
\addlegendentry{$\sim s^3$}
\addplot[red,mark=square*,mark size=2] coordinates {
 (1000,104.9) (2000,104.9) (5000,104.9) (10000,104.9) (20000,104.9) (50000,104.9) (100000,104.9)
};
\addlegendentry{FFTW}
\addplot[cyan,mark=*,mark size=3] coordinates {
 (1000,2.10) (2000,7.10) (5000,41.58) (10000,196.27) (20000,824.53) (50000,1193.44)
};
\addlegendentry{A2R1L}
\addplot[magenta,mark=o,mark size=4] coordinates {
 (1000,4.9449) (2000,13.0930) (5000,32.6423) (10000,76.7702) (20000,187.6492) (50000,456.1334) (100000,836.0194)
};
\addlegendentry{MLFFT}
\addplot[blue,mark=diamond,mark size=4] coordinates {
 (1000,3.6087) (2000,7.9916) (5000,21.0891) (10000,45.2936) (20000,98.1258) (50000,261.6088) (100000,471.9718)
};
\addlegendentry{MLFFT6}
  \end{loglogaxis}
\end{tikzpicture}
}
\caption{Number of samples and runtimes with respect to the sparsity for the reconstruction of random 5-dimensional sparse trigonometric polynomials when considering frequencies within $\Gamma=\hat{G}_{32}^5$, $r=1$.}
\label{figure:numerics:rand_poly:samples_runtimes:r1}
\end{figure}

\begin{Example}\label{numerics:rand_poly:samples_runtimes}
(Number of samples and runtimes with respect to sparsity).
As in Example~\ref{numerics:rand_poly:samples_errors},
we set the expansion $N:=32$, the threshold parameter $\delta:=10^{-12}$ and the number of detection iterations $r:=1$.
Since we want to compare our methods with a full $d$-dimensional FFT, we restrict ourselves to the 5-dimensional case, which already yields a search space $\Gamma=\hat{G}_{32}^5$ of cardinality $|\Gamma|=65^5=1\,160\,290\,625$.
For sparsity $s:=\vert\mathrm{supp}\,\hat{p}\vert\in\{1\,000,2\,000,5\,000,10\,000,20\,000,50\,000\}$,
we determine the runtimes of Algorithm~\ref{algo:sfft_a2r1labs} 
as well as the runtimes of Algorithm~\ref{algo:sfft_mlfft} and~\ref{algo:sfft_mlfft6}.
For Algorithm~\ref{algo:sfft_mlfft} and~\ref{algo:sfft_mlfft6}, we additionally determine the runtimes for sparsity $s=100\,000$.
For each of the considered cases, we subtract the time required for sampling.
We repeat the tests 10 times and capture the maximum of the obtained numbers of samples and runtimes.
The results are depicted in Figure~\ref{figure:numerics:rand_poly:samples_runtimes:r1}.
In Figure~\ref{figure:numerics:rand_poly:samples_runtimes:r1:samples}, we show the maximal number of samples with respect to the sparsity~$s$.
We denote the results of Algorithm~\ref{algo:sfft_a2r1labs} using single rank-1 lattices by ``A2R1L'', of Algorithm~\ref{algo:sfft_mlfft} using multiple rank-1 lattices by ``MLFFT'' and of Algorithm~\ref{algo:sfft_mlfft6} using multiple rank-1 lattices by ``MLFFT6''.
Additionally, we plot the number of samples required by a 5-dimensional FFT on a full grid and denote this by ``5d FFT''.
We observe that the number of samples for Algorithm~\ref{algo:sfft_a2r1labs} using rank-1 lattices increases roughly like $\sim s^{1.5}$ in this example, which is less than the theoretical results of $\mathcal{O}(s^2)$ suggest. For sparsity $s=50\,000$, the number of samples almost equals the cardinality of the search space $|\hat{G}_{32}^5|$, which is the number of samples required by a 5-dimensional FFT.
Moreover, the number of samples is distinctly lower for Algorithm~\ref{algo:sfft_mlfft} using multiple rank-1 lattices and dramatically lower for Algorithm~\ref{algo:sfft_mlfft6}.
For both, Algorithm~\ref{algo:sfft_mlfft} and~\ref{algo:sfft_mlfft6}, the number of samples behaves approximately like $\sim s$ which is slightly better than the theoretical results of $\mathcal{O}(s\log s)$ and $\mathcal{O}(s\log^2 s)$ from Section~\ref{sec:method:complexity}, respectively.
\newline
We also compare the runtimes of the different approaches and visualize the results in Figure~\ref{figure:numerics:rand_poly:samples_runtimes:r1:runtimes}.
Here, we denote the runtime of the 5-dimensional FFT computed using the MATLAB function \verb|fftn|, which is based on FFTW \cite{fftw}, by ``FFTW''.
The observed runtime when using Algorithm~\ref{algo:sfft_a2r1labs} (single rank-1 lattices) is smaller for sparsity $s$ up to 5\,000 and behaves approximately like $\sim s^2\log s$, which is distinctly lower than the theoretical results of $\mathcal{O}(s^3+s^2\log s)$ suggest.
Additionally, when using Algorithm~\ref{algo:sfft_mlfft} (multiple rank-1 lattices, ``MLFFT''), we observe that the runtimes behave approximately like $\sim s\log^2 s$, which corresponds to the theoretical results of $\mathcal{O}(s\log^2 s)$ for fixed $r,d,N$. In particular, these runtimes are distinctly smaller compared to using single rank-1 lattices for higher sparsity~$s$.
Furthermore, Algorithm~\ref{algo:sfft_mlfft} was faster than the 5-dimensional FFT for sparsities $s$ up to $10\,000$ as well as faster than Algorithm~\ref{algo:sfft_a2r1labs} for sparsity $s=5\,000$ and above.
Moreover, the runtimes of Algorithm~\ref{algo:sfft_mlfft6} (multiple rank-1 lattices, ``MLFFT6'') seem to behave slightly better than Algorithm~\ref{algo:sfft_mlfft}, whereas the theoretical results suggest $\mathcal{O}(s\log^3 s)$ arithmetic operations.
\newline
We remark that for higher-dimensional search spaces, the direct application of a $d$-dimensional FFT algorithm is usually not possible in practice due to the exponential growth in the degrees of freedom.
\end{Example}

In the previous examples, we chose the threshold parameter $\delta$ very low compared to the minimal allowed size $|\hat{p}_\boldk|\geq 10^{-6}$ of the Fourier coefficients $\hat{p}_\boldk$.
Next, we consider larger threshold parameters $\delta$ and observe that setting the number of detection iterations $r:=1$ does not always suffice. Still, the observed results are very good.

\begin{Example}\label{ex:sfft:dimincr:per:numerics:success_delta_r}
(Absolute threshold parameter $\delta$ and number of detection iterations $r$).
We set the expansion $N:=32$ and the dimension $d:=5$. In contrast to the previous examples, we use Fourier coefficients $\hat{p}_\boldk:=\mathrm{e}^{2\pi\mathrm{i}\varphi_\boldk}\in\C$ with angles $\varphi_\boldk\in[0,1)$ chosen uniformly at random.
For sparsity $\vert\mathrm{supp}\,\hat{p}\vert:=$ 1\,000 and 1\,000 repetitions with newly chosen random frequencies $\boldk\in\Gamma:=\hat{G}_{32}^5$ and Fourier coefficients $\hat{p}_\boldk$, we investigate the minimal number of correctly detected frequencies and the success rate for correctly detecting all frequencies of Algorithm~\ref{algo:sfft_mlfft}, \ref{algo:sfft_mlfft6} and~\ref{algo:sfft_a2r1labs} for different choices of the absolute threshold parameter $\delta$ and number of detection iterations $r$. The results are depicted in Table~\ref{table:sfft:dimincr:per:numerics:success_delta_r}.
We observe that for increasing number of detection iterations and decreasing threshold parameters $\delta$, the minimal number of correctly detected frequencies and the success rate for correctly detecting all frequencies increase. The three considered algorithms perform similarly well. For instance, setting $r:=3$ and $\delta:=10^{-2}$, all 1\,000 frequencies for all 1\,000 repetitions were successfully detected by each of the three algorithms.
\end{Example}

\begin{table}[htb]
\centering
\begin{tabular}{|r||r|r|r||r|r|r||r|r|r|}
\hline
 & \multicolumn{3}{c||}{Algorithm~\ref{algo:sfft_a2r1labs} using}& \multicolumn{3}{c||}{Algorithm~\ref{algo:sfft_mlfft} using} & \multicolumn{3}{c|}{Algorithm~\ref{algo:sfft_mlfft6} using} \\
 & \multicolumn{3}{c||}{single rank-1 lattices}& \multicolumn{3}{c||}{multiple rank-1 lattices} & \multicolumn{3}{c|}{multiple rank-1 lattices} \\
$\delta$ & 3e-1 & 1e-1 & 1e-2  & 3e-1 & 1e-1 & 1e-2  & 3e-1 & 1e-1 & 1e-2 \\
\hline
$r=1$\rule[0.5em]{0em}{0.55em} & 875 / & 948 / & 987 / & 876 / & 949 / & 993 / & 859 / & 954 / & 987 / \\
      & 0.000 & 0.042 & 0.747 & 0.000 & 0.037 & 0.726 & 0.000 & 0.041 & 0.709 \\
\hline
$r=2$\rule[0.5em]{0em}{0.55em} & 975 / & 984 / & 998 / & 979 / & 996 / & 998 / & 982 / & 984 / & 1000 / \\
      & 0.395 & 0.913 & 0.999 & 0.409 & 0.919 & 0.998 & 0.386 & 0.915 & 1.000 \\
\hline 
$r=3$\rule[0.5em]{0em}{0.55em} & 996 /  & 998 / & 1000 / & 996 / & 998 / & 1000 / & 995 / & 998 / & 1000 / \\
      & 0.912 & 0.996 & 1.000 & 0.931 & 0.997 & 1.000 & 0.907 & 0.997 & 1.000 \\
\hline 
$r=4$\rule[0.5em]{0em}{0.55em} & 998 /  & 1000 / & 1000 / & 998 / & 1000 / & 1000 / & 997 / & 1000 / & 1000 / \\
      & 0.995 & 1.000 & 1.000 & 0.988 & 1.000 & 1.000 & 0.992 & 1.000 & 1.000 \\
\hline
$r=5$\rule[0.5em]{0em}{0.55em} & 998 / & 1000 / & 1000 / & 1000 / & 1000 / & 1000 / &  1000 / & 1000 / & 1000 / \\
      & 0.997 & 1.000 & 1.000 & 1.000 & 1.000 & 1.000 & 1.000 & 1.000 & 1.000 \\
\hline 
\end{tabular}
\caption{Minimal number of correctly detected frequencies (first number) and success rate for correctly detecting all frequencies (second number) for different threshold parameters~$\delta$ and numbers of detection iterations~$r$ when applying Algorithm~\ref{algo:sfft_mlfft}, \ref{algo:sfft_mlfft6}, \ref{algo:sfft_a2r1labs} on 5-dimensional random trigonometric polynomials $p$ of sparsity $\vert\mathrm{supp}\,\hat{p}\vert:=$ 1\,000.}
\label{table:sfft:dimincr:per:numerics:success_delta_r}
\end{table}

\subsection{Random sparse trigonometric polynomial with complex Gaussian noise}
\label{sec:numerics:noise}

In this subsection, we test the robustness to noise. To this end, we build sparse multivariate trigonometric polynomials $p$ with frequencies $\boldk$ randomly chosen from the cube $\hat{G}_N^d:=[-N,N]^d\cap\Z^d$ and corresponding Fourier coefficients $\hat{p}_\boldk:=\mathrm{e}^{2\pi\mathrm{i}\varphi_\boldk}\in\C$ with angles $\varphi_\boldk\in[0,1)$ chosen uniformly at random. We perturb the sampling values of the trigonometric polynomial $p$ by additive complex white Gaussian noise $\eta_j\in\C$ with zero mean and standard deviation $\sigma$,
i.e., we have measurements $\tilde{p}(\boldx_j) = p(\boldx_j) + \eta_j$.
Then, we may approximately compute the signal-to-noise ratio (SNR) in our case by
$$
\mathrm{SNR} \approx \frac{\sum_{j=0}^{M-1} \vert p(\boldx_j)\vert^2 / M}{\sum_{j=0}^{M-1} \vert\eta_j\vert^2 / M}
\approx \frac{\sum_{\boldk\in\mathrm{supp}\,\hat{p}} \vert\hat{p}_\boldk\vert^2}{\sigma^2} = \frac{\vert\mathrm{supp}\,\hat{p}\vert}{\sigma^2}
$$
and we choose $\sigma:=\sqrt{\vert\mathrm{supp}\,\hat{p}\vert}/\sqrt{\mathrm{SNR}}$ for a targeted SNR value.
For our numerical tests in MATLAB, we generate the noise by
$\eta_j := \sigma / \sqrt{2}$\verb| * (randn + 1i*randn)|, $j=0,\ldots,M-1$.
The SNR is often measured using the logarithmic decibel scale (dB), where
$\mathrm{SNR_{dB}} = 10 \log_{10} \mathrm{SNR}$ and $\mathrm{SNR} = 10^{\mathrm{SNR_{dB}}/10}$,
i.e., a linear $\mathrm{SNR}=10^8$ corresponds to a logarithmic $\mathrm{SNR_{dB}}=80\mathrm{dB}$
and $\mathrm{SNR}=1$ corresponds to $\mathrm{SNR_{dB}}=0\mathrm{dB}$.

First, we consider the case of Algorithm~\ref{algo:sfft_a2r1labs} using single reconstructing rank-1 lattices in Example~\ref{ex:numerics:noise:results:A2R1L}. Afterwards, we compare the results with Algorithm~\ref{algo:sfft_mlfft} using multiple rank-1 lattices in Example~\ref{ex:numerics:noise:results:MLFFT4}. In doing so, we observe that one can achieve comparable reconstruction success and errors using Algorithm~\ref{algo:sfft_mlfft} while requiring distinctly less samples compared to using Algorithm~\ref{algo:sfft_a2r1labs}.

\begin{Example}\label{ex:numerics:noise:results:A2R1L}
(see also \cite[Example~3.15]{PoVo14}, reconstruction from noisy sampling values using Algorithm~\ref{algo:sfft_a2r1labs} and reconstructing single rank-1 lattices).
We set the dimension $d:=10$, the expansion $N:=32$, the sparsity $\vert\mathrm{supp}\,\hat{p}\vert:=1\,000$,
and we use the search space $\Gamma:=\hat{G}_{32}^{10}$.
Moreover, we set the sparsity parameter $s:=1\,000$ and the threshold parameter $\delta:=10^{-12}$.
Algorithm~\ref{algo:sfft_a2r1labs} was run setting the parameter $r$ for the number of detection iterations to $r:=1,2,5$
and using the SNR values $\mathrm{SNR_{dB}}:=80,70,\ldots,10,0$ (which corresponds to $\mathrm{SNR}=10^8,10^7,\ldots,10,1$).
For each of these test settings, we repeated the tests $1\,000$ times with new randomly chosen frequencies and Fourier coefficients.
The numerical results are presented in Table~\ref{table:numerics:noise:results:a2r1labs}.
The total number of samples for each of the $1\,000$ repetitions was computed and the maximum of these numbers for each test setting can be found in the column ``\#samples''.
In the column ``min \#freq. correct'', the minimal number of correctly detected frequencies $\vert I^{(1,\ldots,10)}\cap\mathrm{supp}\,\hat{p}\vert$ for the $1\,000$ repetitions is shown,
where $\mathrm{supp}\,\hat{p}$ denotes the set of true (input) frequencies of a trigonometric polynomial $p$ and $I^{(1,\ldots,10)}$ the frequencies returned by the detection algorithm.
The column ``success rate (all freq. correct)'' represents the relative number of the $1\,000$ repetitions where all frequencies were successfully detected,
$I^{(1,\ldots,10)}=\mathrm{supp}\,\hat{p}$.
Moreover, the relative $\ell_2$-error $\Vert(\tilde{\hat{p}}_\boldk)_{\boldk\in I} - (\hat{p}_\boldk)_{\boldk\in I}\Vert_2 / \Vert(\hat{p}_\boldk)_{\boldk\in I}\Vert_2$
of the computed coefficients
$(\tilde{\hat{p}}_\boldk)_{\boldk\in I^{(1,\ldots,10)}}$
was determined for each repetition,
where $I:=\mathrm{supp}\,\hat{p} \, \cup \, I^{(1,\ldots,10)}$ and $\tilde{\hat{p}}_\boldk:=0$ for $\boldk\in I\setminus I^{(1,\ldots,10)}$,
and the column ``rel. $\ell_2$-error'' contains the maximal value of the $1\,000$ repetitions.
We obtain results similar to \cite[Example~3.15]{PoVo14}.
In general, we observe that for decreasing SNR values, the minimal number of correctly detected frequencies and the success rate decrease.
When using $r=1$ detection iteration, there were always some (of the $1\,000$ test runs),
where a few frequencies were incorrect.
When we increased the number of detection iterations $r$, the SNR level at which all frequencies in all of the $1\,000$ test runs were correctly detected also decreased.
For instance for $r=5$ detection iterations, the success rate was at 100 percent including the case $\mathrm{SNR_{dB}}=\mathrm{SNR}=10$.
However, we require about $4$ times of the samples for $r=5$ detection iterations compared to the tests with $r=1$.
\end{Example}

\begin{table}[htb]
\centering
\begin{small}
\begin{tabular}{|S[table-format=2.0]|c|c||S[table-format=9.0]|S[table-format=4.0]|S[table-format=1.3]|c|}
\hline
\multicolumn{1}{|c|}{$\mathrm{SNR_{dB}}$} & noise $\sigma$ & \multicolumn{1}{c||}{\#detect.} & \multicolumn{1}{c|}{max. total} & \multicolumn{1}{c|}{min \#freq.} & \multicolumn{1}{c|}{success rate} & max. rel. \\
 & & \multicolumn{1}{c||}{iter. $r$} & \multicolumn{1}{c|}{\#samples} & \multicolumn{1}{c|}{correct} & \multicolumn{1}{c|}{(all freq. correct)} & $\ell_2$-error \\
\hline
   80\rule[0.5em]{0em}{0.55em} & 3.2e-03 & 1 &    22580805 &       998 &   0.995 & 4.5e-02 \\ 
   70 & 1.0e-02 & 1 &    22783605 &       998 &   0.991 & 5.5e-02 \\ 
   60 & 3.2e-02 & 1 &    22046375 &       998 &   0.973 & 5.5e-02 \\ 
   50 & 1.0e-01 & 1 &    23930465 &       996 &   0.944 & 6.3e-02 \\ 
   40 & 3.2e-01 & 1 &    22263475 &       994 &   0.760 & 7.7e-02 \\ 
   30 & 1.0e+00 & 1 &    22217585 &       992 &   0.443 & 1.0e-01 \\ 
\hline
   80\rule[0.5em]{0em}{0.55em} & 3.2e-03 & 2 &    41573935 &      1000 &   1.000 & 2.5e-06 \\ 
   70 & 1.0e-02 & 2 &    42387345 &      1000 &   1.000 & 7.4e-06 \\ 
   60 & 3.2e-02 & 2 &    42284645 &      1000 &   1.000 & 2.4e-05 \\ 
   50 & 1.0e-01 & 2 &    42055975 &      1000 &   1.000 & 7.2e-05 \\ 
   40 & 3.2e-01 & 2 &    41578745 &      1000 &   1.000 & 2.3e-04 \\ 
   30 & 1.0e+00 & 2 &    41531555 &       998 &   0.991 & 5.5e-02 \\ 
\hline
   50\rule[0.5em]{0em}{0.55em} & 1.0e-01 & 5 &    94793335 &      1000 &   1.000 & 6.3e-05 \\ 
   40 & 3.2e-01 & 5 &    96350345 &      1000 &   1.000 & 1.9e-04 \\ 
   30 & 1.0e+00 & 5 &    96262595 &      1000 &   1.000 & 6.5e-04 \\ 
   20 & 3.2e+00 & 5 &    95140955 &      1000 &   1.000 & 2.0e-03 \\ 
   10 & 1.0e+01 & 5 &    94813745 &      1000 &   1.000 & 6.4e-03 \\ 
    0 & 3.2e+01 & 5 &    93780115 &       999 &   0.800 & 3.8e-02 \\ 
\hline
\end{tabular}
\end{small}
\caption{Results for random sparse trigonometric polynomials
perturbed by additive white Gaussian noise using Algorithm~\ref{algo:sfft_a2r1labs} (reconstructing single rank-1 lattices).}
\label{table:numerics:noise:results:a2r1labs}
\end{table}

\begin{Example}\label{ex:numerics:noise:results:MLFFT4}
Next, we apply Algorithm~\ref{algo:sfft_mlfft} using the identical setting and parameters from Example~\ref{ex:numerics:noise:results:A2R1L}.
We observe similar results as before in Table~\ref{table:numerics:noise:results:MLFFT_1000}.
In general, Algorithm~\ref{algo:sfft_mlfft} requires about 2/3 of the number of samples compared to Algorithm~\ref{algo:sfft_a2r1labs}
and has slightly better success rates. We suspect that this is due to the computation of the projected Fourier coefficients by Algorithm~\ref{algo:mlfft2:ifft_direct} using averaging.
\newline
Additionally, we apply Algorithm~\ref{algo:sfft_mlfft6} and present the results in Table~\ref{table:numerics:noise:results:MLFFT4_1000}. When using $r=1$ detection iteration, there were again always some test runs (of $1\,000$), where not all frequencies could be exactly detected. Moreover, we observe that for some of the test runs, Algorithm~\ref{algo:sfft_mlfft6} with the used parameters was unable to correctly detect any frequency and the minimal number of correctly detected frequencies ``min \#freq. correct''=0. When we increase the number of detection iterations $r$ to 2, the reconstruction success was distinctly higher but worse than the results of Algorithm~\ref{algo:sfft_a2r1labs} in Table~\ref{table:numerics:noise:results:a2r1labs} and of Algorithm~\ref{algo:sfft_mlfft} in Table~\ref{table:numerics:noise:results:MLFFT4_1000}.
Similarly, for $r=5$ detection iterations, we were able to reconstruct better at higher noise levels but still worse compared to Table~\ref{table:numerics:noise:results:a2r1labs} and~\ref{table:numerics:noise:results:MLFFT_1000}.
\newline
Next, we changed the parameter $s_\mathrm{local}$ of Algorithm~\ref{algo:sfft_mlfft6} to the value 1\,200, which means that more frequency candidates are considered and kept in step~\ref{enum:reco:incr_general:I2:non-zero}. Now the minimal number of correctly detected frequencies is distinctly better in many cases but still worse than using Algorithm~\ref{algo:sfft_a2r1labs} or~\ref{algo:sfft_mlfft}.
We suspect that this is also caused by the distinctly smaller number of used sampling values for Algorithm~\ref{algo:sfft_mlfft6}, which is about 1/9 of Algorithm~\ref{algo:sfft_a2r1labs} and 1/6 of Algorithm~\ref{algo:sfft_mlfft}.
\end{Example}

\begin{table}[htb]
\centering
\begin{small}
\begin{tabular}{|S[table-format=2.0]|c|c||S[table-format=9.0]|S[table-format=4.0]|S[table-format=1.3]|c|}
\hline
\multicolumn{1}{|c|}{$\mathrm{SNR_{dB}}$} & noise $\sigma$ & \multicolumn{1}{c||}{\#detect.} & \multicolumn{1}{c|}{max. total} & \multicolumn{1}{c|}{min \#freq.} & \multicolumn{1}{c|}{success rate} & max. rel. \\
 &  & \multicolumn{1}{c||}{iter. $r$} & \multicolumn{1}{c|}{\#samples} & \multicolumn{1}{c|}{correct} & \multicolumn{1}{c|}{(all freq. correct)} & $\ell_2$-error \\
\hline
   80\rule[0.5em]{0em}{0.55em} & 3.2e-03 & 1 &    13032485 &       998 &   0.998 & 4.6e-02 \\ 
   70 & 1.0e-02 & 1 &    13536429 &       996 &   0.990 & 6.5e-02 \\ 
   60 & 3.2e-02 & 1 &    13310277 &       998 &   0.988 & 4.7e-02 \\ 
   50 & 1.0e-01 & 1 &    12911065 &       994 &   0.949 & 7.8e-02 \\ 
   40 & 3.2e-01 & 1 &    13301373 &       996 &   0.857 & 6.7e-02 \\ 
   30 & 1.0e+00 & 1 &    13405817 &       992 &   0.621 & 9.2e-02 \\ 
\hline
   80\rule[0.5em]{0em}{0.55em} & 3.2e-03 & 2 &    25290563 &      1000 &   1.000 & 4.4e-06 \\ 
   70 & 1.0e-02 & 2 &    24883447 &      1000 &   1.000 & 1.4e-05 \\ 
   60 & 3.2e-02 & 2 &    24589823 &      1000 &   1.000 & 4.4e-05 \\ 
   50 & 1.0e-01 & 2 &    25252543 &      1000 &   1.000 & 1.4e-04 \\ 
   40 & 3.2e-01 & 2 &    24937437 &      1000 &   1.000 & 4.4e-04 \\ 
   30 & 1.0e+00 & 2 &    25433681 &       998 &   0.997 & 4.8e-02 \\ 
\hline
   50\rule[0.5em]{0em}{0.55em} & 1.0e-01 & 5 &    63052299 &      1000 &   1.000 & 1.4e-04 \\ 
   40 & 3.2e-01 & 5 &    63618665 &      1000 &   1.000 & 4.3e-04 \\ 
   30 & 1.0e+00 & 5 &    63294069 &      1000 &   1.000 & 1.4e-03 \\ 
   20 & 3.2e+00 & 5 &    65500751 &      1000 &   1.000 & 4.4e-03 \\ 
   10 & 1.0e+01 & 5 &    63537849 &      1000 &   1.000 & 1.5e-02 \\ 
    0 & 3.2e+01 & 5 &    67046015 &       998 &   0.992 & 6.2e-02 \\ 
\hline
\end{tabular}
\end{small}
\caption{Results for random sparse trigonometric polynomials
perturbed by additive white Gaussian noise using Algorithm~\ref{algo:sfft_mlfft} (reconstructing multiple rank-1 lattices).}
\label{table:numerics:noise:results:MLFFT_1000}
\end{table}

\begin{table}[htb]
\centering
\begin{small}
\begin{tabular}{|S[table-format=2.0]|c|c||S[table-format=9.0]|S[table-format=4.0]|S[table-format=1.3]|c|}
\hline
\multicolumn{1}{|c|}{$\mathrm{SNR_{dB}}$} & noise $\sigma$ & \multicolumn{1}{c||}{\#detect.} & \multicolumn{1}{c|}{max. total} & \multicolumn{1}{c|}{min \#freq.} & \multicolumn{1}{c|}{success rate} & max. rel. \\
 &  & \multicolumn{1}{c||}{iter. $r$} & \multicolumn{1}{c|}{\#samples} & \multicolumn{1}{c|}{correct} & \multicolumn{1}{c|}{(all freq. correct)} & $\ell_2$-error \\
\hline
   80\rule[0.5em]{0em}{0.55em} & 3.2e-03 & 1 &     1617311 &         1 &   0.999 & 1.0e+00 \\ 
   70 & 1.0e-02 & 1 &     1616183 &         0 &   0.990 & 1.0e+00 \\ 
   60 & 3.2e-02 & 1 &     1616843 &         0 &   0.967 & 1.0e+00 \\ 
   50 & 1.0e-01 & 1 &     1614591 &         0 &   0.906 & 1.0e+00 \\ 
   40 & 3.2e-01 & 1 &     1617431 &         0 &   0.561 & 1.0e+00 \\ 
   30 & 1.0e+00 & 1 &     1614189 &         0 &   0.000 & 1.0e+00 \\ 
\hline
   80\rule[0.5em]{0em}{0.55em} & 3.2e-03 & 2 &     3063489 &      1000 &   1.000 & 9.3e-06 \\ 
   70 & 1.0e-02 & 2 &     3064963 &      1000 &   1.000 & 3.0e-05 \\ 
   60 & 3.2e-02 & 2 &     3069693 &      1000 &   1.000 & 9.4e-05 \\ 
   50 & 1.0e-01 & 2 &     3130709 &       841 &   0.999 & 4.2e-01 \\ 
   40 & 3.2e-01 & 2 &     3257421 &       570 &   0.934 & 6.8e-01 \\ 
   30 & 1.0e+00 & 2 &     3379167 &       505 &   0.042 & 7.3e-01 \\ 
\hline
   50\rule[0.5em]{0em}{0.55em} & 1.0e-01 & 5 &     7582333 &      1000 &   1.000 & 2.9e-04 \\ 
   40 & 3.2e-01 & 5 &     7570209 &      1000 &   1.000 & 9.1e-04 \\ 
   30 & 1.0e+00 & 5 &     7611301 &       997 &   0.100 & 5.5e-02 \\ 
   20 & 3.2e+00 & 5 &     7855307 &       971 &   0.000 & 1.7e-01 \\ 
   10 & 1.0e+01 & 5 &     9310751 &       834 &   0.000 & 4.1e-01 \\ 
    0 & 3.2e+01 & 5 &    13998685 &         0 &   0.000 & 1.0e+00 \\ 
\hline
\end{tabular}
\end{small}
\caption{Results for random sparse trigonometric polynomials
perturbed by additive white Gaussian noise using Algorithm~\ref{algo:sfft_mlfft6} 
with $s_\mathrm{local}:=1\,000$.}
\label{table:numerics:noise:results:MLFFT4_1000}
\end{table}

\begin{table}[ht]
\centering
\begin{small}
\begin{tabular}{|S[table-format=2.0]|c|c||S[table-format=9.0]|S[table-format=4.0]|S[table-format=1.3]|c|}
\hline
\multicolumn{1}{|c|}{$\mathrm{SNR_{dB}}$} & noise $\sigma$ & \multicolumn{1}{c||}{\#detect.} & \multicolumn{1}{c|}{max. total} & \multicolumn{1}{c|}{min \#freq.} & \multicolumn{1}{c|}{success rate} & max. rel. \\
 & & \multicolumn{1}{c||}{iter. $r$} & \multicolumn{1}{c|}{\#samples} & \multicolumn{1}{c|}{correct} & \multicolumn{1}{c|}{(all freq. correct)} & $\ell_2$-error \\
\hline
   80\rule[0.5em]{0em}{0.55em} & 3.2e-03 & 1 &     1940759 &       995 &   0.999 & 7.7e-02 \\ 
   70 & 1.0e-02 & 1 &     1938367 &       982 &   0.994 & 1.4e-01 \\ 
   60 & 3.2e-02 & 1 &     1942751 &       971 &   0.980 & 1.8e-01 \\ 
   50 & 1.0e-01 & 1 &     1937867 &       969 &   0.955 & 1.8e-01 \\ 
   40 & 3.2e-01 & 1 &     1941323 &       959 &   0.835 & 2.1e-01 \\ 
   30 & 1.0e+00 & 1 &     1939205 &       940 &   0.065 & 2.5e-01 \\ 
\hline
   80\rule[0.5em]{0em}{0.55em} & 3.2e-03 & 2 &     3957003 &      1000 &   1.000 & 8.1e-06 \\ 
   70 & 1.0e-02 & 2 &     3952989 &      1000 &   1.000 & 2.6e-05 \\ 
   60 & 3.2e-02 & 2 &     3972265 &      1000 &   1.000 & 8.2e-05 \\ 
   50 & 1.0e-01 & 2 &     3954775 &      1000 &   1.000 & 2.6e-04 \\ 
   40 & 3.2e-01 & 2 &     3957161 &       998 &   0.984 & 4.5e-02 \\ 
   30 & 1.0e+00 & 2 &     3954715 &       985 &   0.089 & 1.3e-01 \\ 
\hline
   50\rule[0.5em]{0em}{0.55em} & 1.0e-01 & 5 &    10785425 &      1000 &   1.000 & 2.5e-04 \\ 
   40 & 3.2e-01 & 5 &    10768847 &       998 &   0.987 & 4.5e-02 \\ 
   30 & 1.0e+00 & 5 &    10777273 &       994 &   0.105 & 8.4e-02 \\ 
   20 & 3.2e+00 & 5 &    10848099 &       971 &   0.000 & 1.7e-01 \\ 
   10 & 1.0e+01 & 5 &    10804485 &       843 &   0.000 & 4.0e-01 \\ 
    0 & 3.2e+01 & 5 &    17181009 &         0 &   0.000 & 1.0e+00 \\ 
\hline
\end{tabular}
\end{small}
\caption{Results for random sparse trigonometric polynomials
perturbed by additive white Gaussian noise using Algorithm~\ref{algo:sfft_mlfft6} 
with $s_\mathrm{local}:=1\,200$.}
\label{table:numerics:noise:results:MLFFT4_1200}
\end{table}

\FloatBarrier

\newpage
\subsection{Approximation of tensor-product function by trigonometric polynomials}
\label{sec:numerics:fct}

Next, we consider the multivariate periodic test function $f\colon\T^{10}\rightarrow\R$,
\begin{equation} \label{equ:f:10}
f(\boldx):=\prod_{t\in\{1,3,8\}}N_2(x_t) + \prod_{t\in\{2,5,6,10\}}N_4(x_t) + \prod_{t\in\{4,7,9\}}N_6(x_t),
\end{equation}
from \cite[Section~3.3]{PoVo14} with infinitely many non-zero Fourier coefficients,
where $N_m:\T\rightarrow\R$ is the B-Spline of order $m\in\N$,
$$N_m(x) := C_m \sum_{k\in\Z} \operatorname{sinc}\left(\frac{\pi}{m}k\right)^m (-1)^k \,\mathrm{e}^{2\pi\mathrm{i}kx},$$
with a constant $C_m>0$ such that $\Vert N_m \vert L_2(\T)\Vert=1$.
We approximate the function $f$ by multivariate trigonometric polynomials $p$.
For this, we determine a frequency index set
$I=I^{(1,\ldots,10)}\subset\Gamma:=\hat{G}_N^{10}$ and we compute
approximated Fourier coefficients $\tilde{\hat{p}}_\boldk$, $\boldk\in I$, using Algorithm~\ref{algo:sfft_a2r1labs}, \ref{algo:sfft_mlfft} and~\ref{algo:sfft_mlfft6}.
The obtained frequency index sets $I$ should ``consist of'' the union of three lower dimensional manifolds,
a three-dimensional symmetric hyperbolic cross in the dimensions $1,3,8$,
a four-dimensional symmetric hyperbolic cross in the dimensions $2,5,6,10$
and a three-dimensional symmetric hyperbolic cross in the dimensions $4,7,9$.
All tests are performed 10 times and the relative $L_2(\T^{10})$ approximation error
$$
 \Vert f-\tilde{S}_I f\vert L_2(\T^{10})\Vert / \Vert f\vert L_2(\T^{10})\Vert
 =
 \sqrt{\Vert f\vert L_2(\T^{10})\Vert^2 - \sum_{\boldk\in I}\vert\hat{f}_\boldk\vert^2 + \sum_{\boldk\in I}\vert\tilde{\hat{p}}_\boldk-\hat{f}_\boldk\vert^2} / \Vert f\vert L_2(\T^{10})\Vert
$$
is computed,
where $p=\tilde{S}_I f:=\sum_{\boldk\in I} \tilde{\hat{p}}_\boldk \,\mathrm{e}^{2\pi\mathrm{i}\boldk\cdot\circ}$.

Since the number of non-zero Fourier coefficients is not finite and we use a full grid $\hat{G}_N^{10}$ for the search space, we have to limit the number of detected frequencies $\boldk\in I$ and Fourier coefficients~$\tilde{\hat{p}}_\boldk$. For this, we use two different strategies from \cite[Section~3.3]{PoVo14}.
In Section~\ref{sec:numerics:fct:s_sparse}, we set the sparsity parameter~$s$ appropriately.
Moreover, in Section~\ref{sec:numerics:fct:threshold}, we use the threshold parameter~$\delta$ in order to truncate smaller Fourier coefficients.

\subsubsection{s-sparse}
\label{sec:numerics:fct:s_sparse}

First, we use the sparsity input parameter $s\in\N$ to limit the number of detected Fourier coefficients, see also \cite[Section~3.3.1]{PoVo14}.
In doing so, we set the threshold parameter $\delta$ to a very small value, which is $\delta:=10^{-12}$ for the following examples.
In each dimension increment step~$t$, the detected Fourier coefficients and frequencies are truncated and only those frequencies are kept which belong to the $s\ll\vert\Gamma\vert<\infty$ largest Fourier coefficients. Due to possible aliasing effects, we use a larger value $s_\mathrm{local}$ for the intermediate dimension increment steps $t\in\{2,\ldots,d-1\}$ than for the final truncation of the index set $I^{(1,\ldots,d)}$ in dimension increment step $t=d$ in Algorithm~\ref{algo:sfft_a2r1labs}, \ref{algo:sfft_mlfft} and~\ref{algo:sfft_mlfft6}.

\begin{Example} \label{ex:fct:s_sparse:alg2}
($s$-sparse approximate function reconstruction).
We set the expansion $N:=16,32,64$
and we use the full grids $\Gamma:=\hat{G}_N^{10}$ as search space.
Moreover, we set the number of detection iterations $r:=5$.
The used sparsity input parameters $s$ and $s_\mathrm{local}$ as well as the results
of Algorithm~\ref{algo:sfft_a2r1labs} using single rank-1 lattices, of Algorithm~\ref{algo:sfft_mlfft} using multiple rank-1 lattices, and of Algorithm~\ref{algo:sfft_mlfft6} using multiple rank-1 lattices 
are presented in Table~\ref{table:numerics:fct:s_sparse:a2r1l_mlfft4}.
The column ``max. rel. $L_2$-error'' contains the maximum of
the relative $L_2(\T^{10})$ approximation errors $\Vert f-\tilde{S}_I f\vert L_2(\T^{10})\Vert / \Vert f\vert L_2(\T^{10})\Vert$ of the 10 test runs.
The remaining columns have the same meaning as described in Section \ref{sec:numerics:rand_poly}.
We observe that for increasing sparsity parameter, the number of samples increases
while the relative $L_2(\T^{10})$ approximation error decreases.
\newline
When comparing Algorithm~\ref{algo:sfft_mlfft} and~\ref{algo:sfft_a2r1labs}, we observe that both methods yield similar errors. Furthermore, Algorithm~\ref{algo:sfft_mlfft} requires slightly more samples for small sparsity parameters $s$ and about half the number of samples for large ones compared with Algorithm~\ref{algo:sfft_a2r1labs}.
\newline
We observe that the number of samples could be drastically reduced when using Algorithm~\ref{algo:sfft_mlfft6}, 
whereas the relative $L_2(\T^{10})$ approximation errors are comparable for all three algorithms in most cases.
For instance in the case $N=64$ and $s=5\,000$, the maximal total number of samples for Algorithm~\ref{algo:sfft_a2r1labs} (computed over 10 test runs) was about 2.3 billion and for Algorithm~\ref{algo:sfft_mlfft6} only about 159 million samples, which is almost 1/15, achieving comparable relative errors.
Only in 1 test run (out of 10) for expansion $N=64$ and sparsity $s=1\,000$, Algorithm~\ref{algo:sfft_mlfft6} failed to detect some of the frequencies resulting in a relative error of about $4.5\cdot 10^{-1}$ whereas it succeeded for the remaining 9 test runs yielding a relative error of about $1.2\cdot 10^{-2}$, which corresponds to the results of Algorithm~\ref{algo:sfft_a2r1labs}.
\newline
As observed in \cite[Section~3.3.1]{PoVo14}, it is not sufficient to only increase the used sparsity $s$ but the expansion parameter $N$ also needs to be increased.
Using a large expansion $N$ and a small target sparsity $s$ results in the usage of distinctly more samples,
e.g., about 50 million samples for $N=16$ and sparsity $s=1\,000$ compared to about 243 million samples for $N=64$ for Algorithm~\ref{algo:sfft_a2r1labs}
as well as about 8 million samples for $N=16$ and sparsity $s=1\,000$ compared to about 34 million samples for $N=64$ for Algorithm~\ref{algo:sfft_mlfft6}.
This corresponds to the almost linear factor $N$ in the sample complexity, cf.\ Section~\ref{sec:method:complexity}.
\end{Example}

\begin{table}[htb]
\centering
\begin{small}
\begin{tabular}{|r|r||S[table-format=10.0]|r||S[table-format=10.0]|r||S[table-format=9.0]|r|}
\hline
& & \multicolumn{2}{c||}{single rank-1 lattices} & \multicolumn{2}{c||}{multiple rank-1 lattices} & \multicolumn{2}{c|}{multiple rank-1 lattices} \\
& & \multicolumn{2}{c||}{and Algorithm~\ref{algo:sfft_a2r1labs}} & \multicolumn{2}{c||}{and Algorithm~\ref{algo:sfft_mlfft}} & \multicolumn{2}{c|}{and Algorithm~\ref{algo:sfft_mlfft6}} \\
    & & \multicolumn{1}{r|}{max. total} & max. rel. & \multicolumn{1}{r|}{max. total} & max. rel. & \multicolumn{1}{r|}{max. total} & max. rel. \\
$N$ & $s$ & \multicolumn{1}{r|}{\#samples} & $L_2$-error & \multicolumn{1}{r|}{\#samples} & $L_2$-error & \multicolumn{1}{r|}{\#samples} & $L_2$-error \\
\hline
16\rule[0.55em]{0em}{0.55em} &  1\,000 &  50405091 & 1.2e-02 &  62671623 & 1.2e-02 & 8094293 & 1.3e-02 \\
16 &  2\,000 &  128362707 & 4.3e-03 &  127612901 & 3.9e-03 & 15449367 & 5.0e-03 \\
16 &  3\,000 &  222662847 & 3.4e-03 &  176257025 & 3.0e-03 & 22285615 & 4.1e-03 \\
\hline
32\rule[0.55em]{0em}{0.55em} &  1\,000 &  84241365 & 1.2e-02 & 130597513 & 1.2e-02 & 16623363 & 1.3e-02 \\
32 &  2\,000 &  265019105 & 3.4e-03 & 258312555 & 3.4e-03 & 32226797 & 8.8e-03 \\
32 &  3\,000 &  565847035 & 1.7e-03 & 376108443 & 1.6e-03 & 46948295 & 2.2e-03 \\
32 &  4\,000 &  767068055 & 1.4e-03 & 525215599 & 1.2e-03 & 59623261 & 1.6e-03 \\
\hline
64\rule[0.55em]{0em}{0.55em} &  1\,000 &      242940304 & 1.2e-02 &  251626807 & 1.2e-02 & 34236011 & 4.5e-01 \\
64 &  2\,000 &      563767023 & 3.4e-03 &  547034957 & 3.4e-03 & 66914731 & 3.5e-03 \\
64 &  3\,000 &   946215129 & 1.6e-03 &  778060799 & 1.6e-03 & 98926095 & 1.8e-03 \\
64 &  4\,000 &   1266499683 & 9.8e-04 &  1084767689 & 9.7e-04 & 130967911 & 1.1e-03 \\
64 &  5\,000 &    2285094003 & 7.1e-04 &  1318206819 & 6.9e-04 & 158693803 & 9.6e-04 \\
64 & \!\!\! 10\,000 &    5577419619 & 4.3e-04 &  2639274025 & 3.8e-04 & 298886055 & 5.1e-04 \\
\hline
\end{tabular}
\end{small}
\caption{Results for function $f\colon\T^{10}\rightarrow\R$ from \eqref{equ:f:10} when limiting the number of detected frequencies, $s_\mathrm{local}:=2s$.}
\label{table:numerics:fct:s_sparse:a2r1l_mlfft4}
\end{table}

\FloatBarrier

\subsubsection{threshold-based}
\label{sec:numerics:fct:threshold}

Additionally, we consider a variant from \cite[Section~3.3.2]{PoVo14}. 
Now, we use the absolute threshold parameter $\delta\in (0,1)$ for the truncation.
For the search space $\Gamma$, we use the full grid $\hat{G}_N^{10}$ for various expansions $N\in\N$.
Moreover, we set the sparsity input parameter $s:=\vert\Gamma\vert$ and we use $r:=10$ detection iterations, which is twice as high than in Section~\ref{sec:numerics:fct:s_sparse}.
However, this was required for successfully determining the relevant frequency index sets $I$ and achieving low approximation errors.

\begin{Example} \label{ex:fct:threshold:a2r1l_mlfft4}
(Threshold-based approximate function reconstruction).
We apply Algorithm~\ref{algo:sfft_a2r1labs} and~\ref{algo:sfft_mlfft6} on the 10-dimensional test function~\eqref{equ:f:10}.
The parameters and results are shown in Table~\ref{table:numerics:fct:threshold:a2r1l_mlfft4}.
For each parameter combination, we perform the tests 10 times.
For the truncation of the one-dimensional index sets $I^{(t)}$ of frequency candidates for component $t$, $t\in\{1,\ldots,10\}$,
the threshold parameter $\delta:=$``threshold''/10 is used as well as $\delta:=$``threshold'' for all other truncations.
The column ``max. $|I|$'' denotes the maximal number $|I^{(1,\ldots,d)}|$ of Fourier coefficients, which were returned by Algorithm~\ref{algo:sfft_a2r1labs} and~\ref{algo:sfft_mlfft6}, respectively.
The remaining columns have the same meaning as in Example~\ref{ex:fct:s_sparse:alg2}.
We observe that the total numbers of samples are dramatically smaller compared to the results from Section \ref{sec:numerics:fct:s_sparse} while the relative $L_2(\T^{10})$ approximation errors are comparable in most cases for similar numbers $\vert I\vert$ of Fourier coefficients $\tilde{\hat{p}}_\boldk$ used for the approximation $\tilde{S}_I f$ of $f$.
Moreover, the number of samples when using multiple rank-1 lattices and Algorithm~\ref{algo:sfft_mlfft6} is distinctly smaller compared to using single rank-1 lattices and Algorithm~\ref{algo:sfft_a2r1labs}. However, the number $|I^{(1,\ldots,d)}|$ of Fourier coefficients is slightly larger when using Algorithm~\ref{algo:sfft_mlfft6} for ``threshold'' values above $10^{-5}$ and distinctly larger for $10^{-6}$.
\newline
Additionally, we applied Algorithm~\ref{algo:sfft_mlfft} and we observe a similar behavior as in Example~\ref{ex:fct:s_sparse:alg2}.
Algorithm~\ref{algo:sfft_mlfft} and~\ref{algo:sfft_a2r1labs} yield similar errors. Moreover, Algorithm~\ref{algo:sfft_mlfft} requires slightly more samples for large thresholds parameters~$\delta$ and about half the number of samples for small~$\delta$ compared to Algorithm~\ref{algo:sfft_a2r1labs}.
\end{Example}

\begin{table}[htb]
\centering
\begin{small}
\begin{tabular}{|r|r||S[table-format=5.0]|S[table-format=10.0]|r||S[table-format=5.0]|S[table-format=10.0]|r|}
\hline
&  & \multicolumn{3}{c||}{single rank-1 lattices} & \multicolumn{3}{c|}{multiple rank-1 lattices} \\
&  & \multicolumn{3}{c||}{and Algorithm~\ref{algo:sfft_a2r1labs}} & \multicolumn{3}{c|}{and Algorithm~\ref{algo:sfft_mlfft6}} \\
&  & \multicolumn{1}{c|}{max.} & \multicolumn{1}{c|}{max. total} & max. rel. & \multicolumn{1}{c|}{max.} & \multicolumn{1}{c|}{max. total} & max. rel. \\
$N$ & threshold 
& \multicolumn{1}{c|}{$\vert I\vert$} & \multicolumn{1}{c|}{\#samples} & \multicolumn{1}{c||}{$L_2$-error} & \multicolumn{1}{c|}{$|I|$} & \multicolumn{1}{c|}{\#samples} & $L_2$-error \\
\hline
\hline
64\rule[0.55em]{0em}{0.55em} &   1.0e-02 &   531 &      388370 & 4.2e-01 &  573 &      299551 & 1.9e-01 \\
64 &   1.0e-03 &  1157 &     4089813 & 1.1e-02 &  1327 &     2172161 & 1.1e-02 \\
64 &   1.0e-04 &  3075 &    22232568 & 1.6e-03 &  3847 &     8924419 & 2.1e-03 \\
64 &   1.0e-05 &  8041 &   163970316 & 4.6e-04 &  12499 &    42668463 & 5.1e-04 \\
64 &   1.0e-06 & 22019 &  1087478506 & 4.1e-04 &  56375 &   249240643 & 4.2e-04 \\
\hline
128\rule[0.55em]{0em}{0.55em} & 1.0e-02 &  531 &      413506 & 1.3e-01 &      567 &      328595 & 6.1e-02 \\
128 & 1.0e-03 &  1163 &     3196264 & 1.0e-02 &      1357 &     1865193 & 3.9e-02 \\
128 & 1.0e-04 &  3099 &    32334421 & 1.6e-03 &     3733 &    13946741 & 1.6e-03 \\
128 & 1.0e-05 &  8491 &   250367814 & 3.0e-04 &     11039 &    52172733 & 2.6e-03 \\
128 & 1.0e-06 & 23609 &  1679136524 & 1.5e-04 &   50571 &   316908069 & 1.6e-04 \\
\hline
256\rule[0.55em]{0em}{0.55em} &   1.0e-02 &   531 &      456113 & 2.0e-01 &      571 &      352607 & 6.9e-02 \\
256 &   1.0e-03 &  1155 &     3185671 & 1.1e-02 &    1355 &     1890573 & 1.1e-02 \\
256 &   1.0e-04 &  3085 &    51682418 & 1.6e-03 &    3653 &    18860827 & 1.6e-03 \\
256 &   1.0e-05 &  8553 &   422552990 & 2.9e-04 &    10863 &    77359193 & 2.9e-04 \\
256 &   1.0e-06 & 25563 &  3366629304 & 6.3e-05 &  40331 &   373569439 & 1.0e-04 \\
\hline
\end{tabular}
\end{small}
\caption{Results for approximation of function $f\colon\T^{10}\rightarrow\R$ from~\eqref{equ:f:10} when truncating the detected frequencies and Fourier coefficients by threshold.}
\label{table:numerics:fct:threshold:a2r1l_mlfft4}
\end{table}

\FloatBarrier

\section{Conclusion}
\label{sec:conc}

In this work, we presented fast, efficient, reliable algorithms
for high-dimensional sparse FFTs based on dimension-incremental projections.
We applied the new algorithms to the reconstruction and approximation of high-dimensional functions,
assuming that the frequencies of important Fourier coefficients are unknown.

One main ingredient was the utilization of multiple rank-1 lattices, which are
newly developed sampling schemes allowing for highly efficient construction methods and corresponding FFT algorithms
for known support~$\mathrm{supp}\,\hat{p}$ in frequency domain.

Another important aspect is the modification of algorithm parameters that allow for
theoretical estimates.
In particular, a detailed analysis of specific realizations of the dimension-incremental approach
was performed for the first time with respect to success probability.
Moreover, we estimated the sample and arithmetic complexities
for the case of successful detection as well as for the worst case.

Various numerical tests complement the theoretical considerations and indicate
that one may obtain smaller bounds on the complexities, especially with respect to
the sparsity of the signal under consideration. In Section~\ref{sec:theo:discussion}
we discussed one improvement idea and point out the emerging
open questions.

\subsection*{Acknowledgements}
LK and DP gratefully acknowledge the funding by the Deutsche Forschungsgemeinschaft (DFG, German Research Foundation) -- 380648269. \newline
Additionally, DP and TV gratefully acknowledge the funding by the European Union and the Free State of Saxony (EFRE/ESF NBest-SF).

\begin{small}
\setlength{\bibsep}{0.45\baselineskip}

\end{small}

\end{document}